\documentclass[12pt]{amsart}
\usepackage{amscd,amssymb,verbatim}
\theoremstyle{plain}

\usepackage{color}

\numberwithin{equation}{section}

\newcommand{\N}{\mathbb{N}}
\newcommand{\R}{\mathbb{R}}
\newcommand{\Q}{\mathbb{Q}}

\newcommand{\e}{\varepsilon}
\newcommand{\mc}{\mathcal}

\newcommand{\supp}{\textrm{supp}}
\newcommand{\sign}{\textrm{sign}}

\newtheorem{theorem}{Theorem}[section]
\newtheorem{lemma}[theorem]{Lemma}
\newtheorem{remark}[theorem]{Remark}

\newtheorem{definition}[theorem]{Definition}
\newtheorem{proposition}[theorem]{Proposition}
\newtheorem{corollary}[theorem]{Corollary}
\newtheorem{notation}[theorem]{Notation}

\newtheorem{example}{Example}

\begin{document}

\title[The Extreme Points of the Unit Ball of $J$ and its dual spaces]{The Extreme
Points of the Unit Ball of the James space $J$ and its dual spaces}

\author{Spiros A. Argyros}
\address{National Technical University of Athens, Faculty of Applied Sciences, 
Department of Mathematics, Zografou Campus, 157 80, Athens, Greece.}
\email{sargyros@math.ntua.gr}

\author{Manuel Gonz\'alez}
\address{Department of Mathematics, Faculty of Sciences, University of Cantabria, 
Avenida de los Castros s/n, 39071 Santander, Spain.}
\email{manuel.gonzalez@unican.es}



\begin{abstract} 
We  provide a new proof of S. Bellenot's characterization of the extreme points of the unit ball $B_J$ 
of the James quasi-reflexive space $J$. We also provide an explicit description of the norm of $J^{**}$ 
which yields an analogous characterization for the extreme points of $B_{J^{**}}$. 
In the last part of the paper we describe the  set of all extreme points of $B_{J^*}$ and its norm closure.
It is remarkable that the descriptions of the extreme points of $B_J$ and $B_{J^*}$ are closely connected.
\end{abstract}

\maketitle


\section{Introduction}\label{sect:intro}
The James space, constructed by R.C. James (\cite{J1}) more than seventy years ago, showed the existence 
of a non-reflexive Banach space isomorphic to its second dual, and a quasi-reflexive, non-reflexive 
Banach space. This space played a central role in the development of Banach space theory and was extensively 
studied by many authors. We mention James' paper \cite{J2} and also the works of A. Andrew (\cite{A1}, \cite{A2}), 
S. Bellenot (\cite{B}, \cite{B2}, \cite{B3}) and P.G. Casazza et al.\ (\cite{C}, \cite{CLL}).
A comprehensive presentation of the properties of James' space and its dual space can be found in the monograph 
(\cite{F}) by H. Fetter and B. Gamboa de Buen.

In (\cite{J1}), James' space was constructed using the squared variation norm (Equation \ref{s-norm}) which makes 
the unit vector basis $(e_n)_{n\in \N}$ a shrinking (Schau\-der) basis. We denote this space by $J_s$. 
In this paper, we consider a space $J$ isometric to $J_s$ for which $(e_n)_{n\in \N}$ is a boundedly complete 
basis. It is defined as follows:

A sequence of real numbers $x = (x(n))_{n\in \N}$ is in $J$ when $\|x\|_J<\infty$, where
\begin{equation}\label{eq:Jnorm}
\| x \|_{J} = \sup \biggl\{\biggl( \sum_{i=1}^{n} \big| \sum _{k \in I_{i} }x(k)\big|^2 \biggr)^{1/2} :\; \{ I_i \}_{i=1}^n \text{ disjoint intervals of }
\N \biggr \}.
\end{equation}

This norm appeared later in the literature, and R.C. James (\cite{J2}) attributed its definition to J. Lindenstrauss (\cite{L1}).

Since $\|x\|_2 = \left(\sum_{k=1}^\infty |x(k)|^2 \right)^{1/2}\leq \|x\|_J$, the space $J$ is a dense linear subspace of $\ell_2$.
Moreover, $(e_n)$ is \emph{spreading} for this norm, in the sense that it is isometrically equivalent
to each of its subsequences.

In this paper, we characterize the extreme points of the unit balls of $J$, $J^*$ and $J^{**}$.

In Section \ref{sect:prelim}, we present some properties of $J$ and $J_s$, and we introduce some notation.
Section \ref{sect:ext-B_J} is devoted to  the following characterization of the set of extreme points of the unit ball of $J$, denoted $Ext(B_J)$, first proved by S. Bellenot (\cite{B2}):
A vector $x \in J$ is in $Ext(B_J)$ if and only if $\|x\|_J= \|x\|_2=1$ (Theorem \ref{T1}).

As a consequence of this characterization,  $Ext(B_J)$ is norm-closed (Corollary \ref{C31}).
Theorem \ref{T1} is not explicitly stated in \cite{B2}, where the set of extreme points of the unit ball of $J_s$ was characterized (see Corollary \ref{cor3.2}).
For doing this, Bellenot developed a maximum path algorithm for graphs using a language vaguely related to that of Banach space theory. Here we follow a different approach using the concept of \emph{$x$-norming partition} defined below.

Given $x=(x(n))\in J$, the \emph{support of $x$} is the set $\supp(x)= \{n\in\N : x(n)\neq 0\}$.
%
%
We consider families $\mc I=\{I_i\}_{i\in F}$ of intervals of $\N$, where the set of indices $F$ is $\{1, \ldots,k\}$ for some $k\in\N$ or $F=\N$, each $I_i$ is non-empty, and $\max I_i < \min I_{i+1}$ when $i+1\in F$. So the intervals are pairwise disjoint. We denote
$$
\|x\|_{\mc I}= \bigl(\sum_{i\in F} \bigl|\sum_{k\in I_i} x(k) \bigr|^2 \bigr)^{1/2}.
$$

Given $x\in J$, we say that $\mc I$ is an  \emph{$x$-norming family} if $\|x\|_{\mc I}= \|x\|_J$.
Clearly every  $\mc I=\{I_i\}_{i\in F}$ \emph{$x$-norming family} satisfies $\supp(x) \subset \cup_{i\in F} I_i$.
Next for $\emptyset \ne L\subset \N $ we set $\sup L = \max L$ if $L$ is finite and $\sup L =\infty$ otherwise.

\begin{definition}\label{x-norm-part}
Let $x\in J$. An \emph{$x$-norming partition} is an $x$-norming family $\mc {I}=\{I_i\}_{i\in F}$ of intervals of $\N$ such that for every $i <\sup F$ we have $\{\min I_i, \max I_i\}\subset \supp(x)$, and for $i = \sup F < \infty$, $\min I_i \in \supp(x) $ and $\sup I_i = \sup \supp(x)$.
\end{definition}

Every $x$-norming partition defines a partition of  $\supp(x)$.
An $x$-norming family $\mc J=\{J_j \}_{j\in G}$ has associated an $x$-norming partition obtained as follows:
denoting $J_j^s= J_j\cap \supp(x)$, we delete those $J_j$ with $J_j^s =\emptyset$, replace the remaining ones by
$I_j=[\min J_j^s, \max J_j^s]$ when $J_j^s$ is bounded and by $I_j=[\min J_j^s, \infty)$ otherwise, and re-enumerate
them to obtain $\{I_i\}_{i\in F}$.
We show that every $x\in J$ admits an $x$-norming partition (Corollary \ref{C0}), although it can admit more than one (Section \ref{sect:multiple-part}).

When we discovered our first proof of Theorem \ref{T1}, we were unaware that Bellenot had already proved it. Our proof was based on some properties of the $x$-norming families (Proposition \ref{prop3.1}). Bellenot's concept of "maximal cost paths" helped us to reveal some deeper regularity properties.

The key result (Corollary \ref{cor3.3}) shows that given two $x$-norming partitions $\mathcal{I}$ and $\mathcal{L}$, if $I\in \mathcal{I}$ and $L\in \mathcal{L}$, then  $I \subset L$ or $L\subset I$ or $I\cap L= \emptyset$.
As a consequence, there exists another $x$-norming partition $\mc J$ that refines both $\mathcal{I}$ and $\mathcal{L}$ (Corollary \ref{cor3.4}), and an $x$-norming partition $\mc M$ that refines any other $x$-norming partition (Proposition \ref{prop3.5.1}). We call $\mathcal{M}$ the \emph{finest $x$-norming partition}. It corresponds to Bellenot’s maximal cost path and its existence is the main result of [5]. Our proof of Theorem 3.13 is based on these regularity properties.

Given $x\in J$ with $\|x\|_J=1$, every $x$-norming partition $\{I_i\}_{i\in F}$ generates a point in $Ext(B_J)$ as follows. If $b_i= \sum_{n\in I_i}x(n)$ for $i\in F$, then $\sum_{i\in F}b_i e_{k_i}\in Ext(B_J)$ for every strictly increasing sequence $(k_i)_{i \in F}$ in $\N$  (Corollary \ref{C32}).
We also describe $Ext(B_J)$ without explicitly using the norm of $J$. Given $x\in \ell_2$ with $\|x\|_2=1$, $x\in Ext(B_J)$ if and only if $ \bigl |\sum_{n\in I} x(n)\bigr|^2 \leq \sum_{n\in I} |x(n)|^2 $ for every interval $I$ of $\N$ (Corollary \ref{C33}).

In Section \ref{sect:J**} we consider the set $\mathcal{D}\subset B_{J^*}$ (see Equation \ref{eq:D-set}), we define $\mathcal{D}_1$ as the $w^*$-closure of $\mathcal{D}$ in $J^*$, and we show that $\mathcal{D}_1$ norms $J^{**}$. This fact allows us to give an explicit description of the norm of $J^{**}$, which is a natural extension of the norm of $J$. Denoting $e_\omega= w^*$-$\lim_n e_n\in J^{**}$, the space $J^{**}$ is generated by the transfinite sequence $(e_p)_{p\in \omega +1}$ with the following norm:
$$
\|x^{**} \|_{**} = \sup\biggl\{ \bigl( \sum_{i\in F } |\sum_{n\in I_{i}} x^{**}(n)|^{2} \bigr)^{1/2} : \{I_{i} \}_{i\in F} \ \text{ disjoint intervals of} \ \omega+1 \biggr \}.
$$

In a similar way, we can define the norm of the transfinite even dual spaces of $J$.
Moreover, the points of $Ext(B_{J^{**}})$ are described  similarly to those of $Ext(B_J)$.

Section \ref{sect:J*} is devoted to the study of $Ext(B_{J^*})$. Note that for a non empty interval $I$ of $\N$,
$I^*(x) = \sum_{n\in I}x(n)$ defines $I^*\in J^*$ with $\| I^* \|_* =1$.
Throughout this paper $I^*$ will denote a functional as before with $I$ a non-empty interval of $\N$.

We give a description of the set $\mathcal{D}_1$ (Lemma \ref{lem5.1}), which is a $w^*$-closed subset of $B_{J^*}$
that norms $J$ (\cite{J2}). Hence $Ext(B_{J^*}) \subset \mathcal{D}_1$ by Theorem \ref{thm:H}.

As we have mentioned, the characterization of $Ext(B_{J^*})$ is closely related to that of $Ext(B_J)$. Proposition \ref{prop5.1} and Theorem \ref{thm5.1} indicate this fact.
\smallskip

Proposition \ref{prop5.1}: 
Let $x^*= \sum_{i\in F }\alpha_i I_i^* \in \mathcal{D}_1$. Then $\|x^*\|_*= 1$ if and only if $\sum_{i\in F}\alpha_i e_i$
is an extreme point of $B_J$.
\smallskip

Therefore every $x^*= \sum_{i\in F} \alpha_{i }I_{i}^* \in Ext(B_{J^*})$ must satisfy the above condition.
Furthermore, we show that for such an $x^*$, the set $\cup_{i\in F} I_i$ is an interval of $\N$ (Proposition \ref{prop5.2}).
These two requirements characterize the extreme points of $B_{J^*}$:
\smallskip

Theorem \ref{thm5.1}: A vector $x^* \in J^*$ is in $Ext(B_{J^*})$ if and only if it satisfies

(i)  $x^* = \sum_{i\in F} \alpha_i I_i^*$ with $ \alpha_i\ne 0$ for every $i\in F$, and

(ii)  $x = \sum_{i\in F} \alpha_i e_i \in Ext(B_J)$ and  $\cup_{i\in F} I_i$ is an interval of $\N$.
\smallskip

In contrast to $Ext(B_J)$, $Ext(B_{J^*})$ is not norm closed. In Proposition \ref{prop5.8} we characterize the elements
of the norm-closure of $Ext(B_{J^*})$.

In Section \ref{sect:multiple-part} we give examples of vectors $x\in J$ such that the cardinal of the set of $x$-norming
partitions takes all possible values: any $k\in\N$, $\aleph_0$ or $\mathfrak{c}$ (Proposition \ref{prop00}).

\section{Preliminaries}\label{sect:prelim}

\subsection {Some properties of $J$ and $J^*$}
James space \cite{J1} can be described using two norms. The first norm is the original one in \cite{J1}. It is called
the \emph{square variational norm,} and it is defined as follows.
For $x\in c_{00}$ we set
\begin{equation}\label{s-norm}
\|x\|_{s} = \frac{1}{\sqrt{2}}\max \biggl\{\bigl (\sum_{k=1}^{m} \bigl|x(n_{k})- x(n_{k-1})\bigr|^2 \bigr)^{1/2}
: 1\leq n_1<\cdots< n_m \biggr \},
\end{equation}
where we assume $n_0=0$ and $x(0)=0$.

We denote by $J_s$ the completion of $(c_{00}, \| \cdot \|_s)$, where $c_{00}$ is the space of finitely non-zero real sequences.
The unit vector basis $(e_n)_{n\in \N}$ is a shrinking basis for $J_s$.

The second norm $\|\cdot\|_J$ (Equation \ref{eq:Jnorm}) appeared for the first time in  (\cite{J2}). James attributed its
definition to J. Lindenstrauss (\cite{L1}). The unit vector basis $(e_n)_{n\in \N}$ is a boundedly complete basis for $J$.

A vector $x = \sum_n a_n e_n \in J$ will often be denoted as a scalar sequence $x \colon \N \to \mathbb{R}$, where $x(n) = a_n $ for all $n \in \N$.

The spaces $J$ and $J_s$ are isometric. Indeed, it is easy to check that the linear map $T:J \to J_s$ defined by $T(x)= (\sum_{k=n}^{\infty}x(k))_n$ is an onto isometry. The inverse $T^{-1}$ maps $(a_n)_n \in J_s$ to $(a_n- a_{n+1})_n\in J$.

James \cite{J2} showed that $J^*$ is not isomorphic to any subspace of $J$, and Andrew \cite{A1} showed that no non-reflexive subspace of $J$ is isomorphic to a subspace of $J^*$.
We also mention two important results concerning the local structure of $J$ and $J^*$: James \cite{J2} proved that $\ell_\infty$ is finitely representable in $J$, and G. Pisier \cite{P} showed that $J^*$ has cotype $2$, hence $\ell_\infty$ is not finitely representable in $J^*$.

\subsection{Notation}
For a Banach space $X$, $B_X$ and $S_X$ are the closed unit ball and the unit sphere of $X$.
The \emph{support of $x\in \R^\N$} is the set $\supp(x) =\{ n\in \N: x(n) \neq 0 \}$, and given $A\subset \N$, we denote by $x|_A$ the sequence that coincides with $x$ for $n\in A$ and it is $0$ otherwise. We also denote by $\# C$ the cardinal of a set $C$.

If $\mathcal{I}=\{I_i\}_{i\in F}$ is a family of disjoint intervals of $\N$, we assume that $\max I_i < \min I_{i+1}$ for $i+1 \in F$, where $F$ is $\{1,\ldots,k\}$ or $\N$.
If $x \in J$ and $\mathcal{I} = (I_i)_{i\in F}$ is a family of disjoint intervals of $\N$, we denote
$$
\|x\|_{\mathcal{I}} =\bigl(\sum_{i\in F} \bigl| \sum_{n\in I_i} x(n)\bigr|^2 \bigr )^{1/2}.
$$

We write $(J^*, \|\cdot\| _*)$ for the dual of $J$. Since $(e_n)$ is a boundedly complete basis for $J$, the sequence of biorthogonal functionals $(e_{n}^*)_n\subset J^*$ is a shrinking basis for its closed linear span in $J^*$.
We denote this subspace of $J^*$ by $J_*$ because it is a predual of $J$.


\section{A characterization of the extreme points of $B_J$}\label{sect:ext-B_J}

\subsection{The existence of $x$-norming partitions}
In this subsection we prove that every $x\in J$ admits at least one $x$-norming partition. If $\supp(x)$ is finite, this is clear. If $L=\supp(x)$ is an infinite subset of $\N$ then every $x$-norming partition $\mc{I}$ defines a partition of $L$. For simplicity, we shall work for the case $L=\N$. Since $(e_n)$ is spreading, the general case follows the same arguments.

We denote by $\mathcal{P}_I(\N)$ the set of all partitions of $\N$ into intervals. Observe that
$\|x\|_J= \sup\{ \|x\|_{\mathcal{I}} : \mathcal{I} \in \mathcal{P}_I(\N)\}$ for every $x\in J$.
Our first task is to show that every $x \in J$  with $\supp(x) = \N$ admits an $x$-norming partition.

\begin{definition}\label{def:set}
We set $m_0=0$, and for every subset $A$ of $\N$ we define a partition $\mathcal{P}_A \in \mathcal{P}_I(\N)$ as follows:
\begin{enumerate}
\item If $A = \emptyset$, then $\mathcal{P}_A=\{\N\}$.
\item If $A=\{m_1<\cdots<m_k\}$ is non-empty and finite, then
$\mathcal{P}_A=\{I_i\}_{i=1}^{k+1}$ with  $I_i=(m_{i-1},m_i]$ for $i\leq k$ and $I_{k+1}=(m_k,\infty)$.
\item If $A=\{m_1< m_2<\cdots\}$ is infinite, then $\mathcal{P}_A= \{I_i\}_{i \in\N}$ with $I_i=(m_{i-1},m_i]$.
\end{enumerate}
Moreover, we denote $\|x\|_A= \|x\|_{\mc{P}_A}$ for $x\in J$.
\end{definition}

Observe that $\|x\|_J =\sup_{A\subset\N}\|x\|_A$ for every $x \in J$. In our next result, by identifying sets with
their indicator functions, we endow the powerset of $\N$ with the topology of pointwise convergence.

\begin{proposition} \label{P0}
Let $x\in J$ with $\supp(x) = \N$ and  $(A_n)$ be a sequence of non-empty, finite subsets of $\N$ pointwise
convergent to $A\subset \N$. Then $\lim_{n\to\infty} \|x\|_{A_n} = \|x\|_A$ for every $x \in J$.
\end{proposition}
\begin{proof}
Let $x\in J$. Let $P_m ,Q_m$ denote the natural projections onto the closed subspaces generated by
$\{e_i :i\le m\}$, $\{e_i : i>m\}$ respectively. Then $\|x\|_A \le \|x\|_J$, hence
$\lim_{m\to \infty} \|Q_mx\|_A=0$ for every $A\subset \N$.
We consider three cases:
\medskip

$A=\emptyset$. In this case $\lim_{k\to \infty} \min A_k=\infty$. Hence
$$
\lim_{k\to \infty}\|x\|_{A_k}= \big| \sum_{n=1}^\infty x(n)\big| =\|x\|_A.
$$

\emph{$A$ finite nonempty.}
Denoting $m_0 = \max A$, we may assume that $A= A_k\cap [1, m_0]$ for each $k$.
Then $\|x\|_A^2= \|P_{m_0}(x)\|_A^2 + \| Q_{m_0}(x)\|_A^2$ and
$$
\|x\|_{A_k}^2= \|P_{m_0}(x)\|_{A}^2+ \|Q_{m_0} (x)\|_{A_k}^2= \| P_{m_0}(x)\|_{A}^2+ \|Q_{m_0}(x)\|_{A_k\setminus A}^2.
$$

Since $\lim_{k\to \infty} \min A_k \setminus A= \infty$,  $\lim_{k\to \infty}\|Q_{m_0}(x) \|_{A_k\setminus A}= \|Q_{m_0}(x)\|_A$.
Hence $\lim_{k\to \infty}\|x\|_{A_k}= \|x\|_A$.
\medskip

\emph{$A$ infinite.} For each $m\in\N$ there exists $k_m$ such that $A_k\cap[1,m]= A\cap[1,m]$ for $k>k_m$.
Hence $\|P_m (x)\|_{A_k} =\|P_m (x) \|_A$ for $k > k_m$.

Given $\e> 0$ we choose $m\in \N$ such that  $\|Q_m(x)\|_J< \e$.
Then, for $k> k_m$, $\|x\|_A< \|P_m(x)\|_{A_k} +\e \le\|x\|_{A_k}+ \e$ and $\|x\|_{A_k}< \|P_m(x)\|_A+\e \le \|x\|_F +\e$,
which yields  $\lim_{k\to \infty}\|x\|_{A_k}= \|x\|_A$.
\end{proof}

As mentioned above, the same proof works for $x\in J$ with $L =\supp(x)$ an infinite subset of $\N$, by considering partitions of $L$ instead of $\N$. 

\begin{corollary} \label{C0}
For every nonzero $x \in J$ there exists $A \subset \N$ such that $\|x\|_A = \|x\|_J$.
Hence there exists an $x$-norming partition.
\end{corollary}
\begin{proof}
By the definition of $\|\cdot\|_J$ there exists a sequence $(A_n)$ of finite, non-empty, subsets of $\N$  with $A_n \subset \supp(x) $ such that $\lim_n \|x\|_{A_n} = \|x\|_J$. By compactness, there is no loss of generality in assuming that $(A_n)$ is pointwise convergent to some $A\subset \N$. By Proposition \ref{P0}, $\|x\|_A = \|x\|_J$. So $\mc P_A$ defines an $x$-norming partition.
\end{proof}

\begin{remark}\label{rem3.1} 
Let $x\in J$ with $L=\supp(x)$ an infinite set. Then the function $\Phi$ mapping an $x$-norming partition $\mc{P}_A$ to the characteristic function of the set $A$ in $\{0,1\}^L$ is injective, and Proposition \ref{P0} yields that its range is a closed subset of $\{0,1 \}^L$. 
In Section \ref{sect:multiple-part} we will see  vectors $x\in J$ for which the set of $x$-norming partitions is uncountable. The previous argument implies that the cardinal of this set is $\mathfrak c$. 
\end{remark}

\subsection{Properties of $x$-norming families}

Recall that the set of indices $F$ is $\{1, \ldots, k\}$ for some $k\in\N$, or $\N$.

\begin{proposition}\label{prop3.1}
Let $x \in J$ and let $\{I_i\}_{i\in F}$ be an $x$-norming family. Then the following assertions hold:
\begin{enumerate}
\item $\supp(x) \subset \cup_{i\in F}I_i$.
\item  For every subinterval $G$ of $F$, we have 
$$
\sum_{i\in G} |\sum_{n \in I_{i}}x(n)|^2 = \|x|_{\cup_{i\in G}I_{i}} \|_J^2.
$$
\item  Given $i\in F$ and $n\in I_i \cap \supp(x)$, the scalars
$$
\sum_{m \in I_{i} } x(m),\,  \sum_{m \in I_i\, m \leq n} x(m)\, \text{ and } \sum_{m \in I_i, \, m \geq n} x(m)
$$
are all non-zero and have the same sign.

\item  For every $i \in F $,
$\sum_{m \in I_{i} } x(m)$ and $\sum_{m \in I_{i+1} } x(m)$ have alternate signs.

\item Let $i\in F$. If $n_1 < n_2$ are consecutive in $\supp(x)$ and $x(n_1)x(n_2)> 0$, then $\{n_1, n_2\} \subset I_i$ or $\{n_1,n_2\} \cap I_i = \emptyset$.
\end{enumerate}
\end{proposition}
\begin{proof}
(1) It is immediate.

(2) If not, then Corollary \ref{C0} allows us to choose a partition $\{ L_m\}$ of $\cup_{i\in G}I_i$ that is $(x|_{\cup_{i\in G}I_{i}})$-norming. This shows that $\{I_i\}_{i\in F}$ is not $x$-norming.
\smallskip

(3) If $\sum_{m \in I_i, \, m \leq n} x(m) =0$, then $| \sum_{m \in I_i, \, m > n} x(m)| =  \| x|_{I_{i}} \|_J  $, a contradiction, since
$$
x(n) \ne 0 \text{ and } \| x|_{I_{i}} \|_J^2  \geq \bigl |\sum_{m \in I_i, \, m > n} x(m) \bigr |^2 + |x(n)|^2.
$$
Similar arguments show that all three scalars have the same sign.

(4) Observe that $\sum_{m\in I_i} x(m)\neq 0$ for all $i\in F$. So, if the conclusion fails,
$$
\bigl | \sum_{m \in  I_{i}\cup I_{i+1} } x(m) \bigr |^2 > \bigl | \sum_{m \in I_{i} } x(m) \bigr |^2 +  \bigl |\sum_{m \in I_{i+1} } x(m) \bigr |^2
$$
which yields that $\{I_i\}_{i\in F}$ is not $x$-norming.

(5) Otherwise, $n_1= \max I_i\cap \supp(x)$ and $n_2= \min I_{i+1}\cap \supp(x)$ for some $i\in F$, and properties (3) and (4) would produce a contradiction.
\end{proof}

\begin{lemma}\label{l1}
Let $\rho$, $\gamma$ and $\varepsilon$ be scalars with $\varepsilon \ne 0$. Assume that $\rho + \gamma^2 \leq (\gamma + \varepsilon)^2$.
Then
$\rho + (\gamma + \delta)^2 <  (\gamma + \delta + \varepsilon)^2$ for every scalar $\delta $ such that $\delta \varepsilon >0$.
\end{lemma}
\begin{proof}
Developing the terms, we easily get $(\gamma + \delta + \varepsilon)^2-\rho - (\gamma + \delta)^2\geq 2\e\delta$.
\end{proof}

\begin{proposition}\label{p1}
Let $0\neq x = \sum_{n=1}^\infty \lambda_ne_n \in J$, and let $\mathcal{I}=\{ I_i\}_{i\in F}$ and $\mathcal{L}= \{ L_j\}_{j\in G}$ be $x$-norming partitions. Assume that there exist $i\in F$ and $j\in G$ such that
$$
I_i \subset L_j, \text{ and } \min I_i = \min L_j.
$$
Then for every $i' \in F$, $i'\neq i$, either $I_{i'} \subset L_j$, or $I_{i'} \cap L_j = \emptyset$.

Therefore, there exists a subinterval $F_j$ of $F$ such that
\[ L_j \cap \supp(x) = \cup_{i\in F_j}[ I_i \cap \supp(x)].  \]
\end{proposition}
\begin{proof}
Assume that the conclusion fails. We prove the case  $j=1$, since the general case is proved with exactly the same arguments.

(i) Since $j= 1$, we conclude that $i = 1$.

(ii) Let $I_k$ be the unique element of $\{I_i\}_{i \in F}$ that is separated by $L_1$. Since
    \[ L_1 \cap \supp(x) =\bigl[ (\cup_{m < k} I_m)  \cup (L_1 \cap I_k )\bigr] \cap \supp(x), \]
we obtain
\begin{equation}\label{eq:scalars}
    \sum_{n \in L_1} \lambda_n= \sum_{p < k} \sum_{n \in I_p} \lambda_n +   \sum_{n \in L_1 \cap I_k} \lambda_n.
\end{equation}

(iii) Equation \ref{eq:scalars} combined with part (3) of Proposition \ref{prop3.1} yield that the scalars
$$
\sum_{m < k} \sum_{n \in I_m} \lambda_n,  \,\, \sum_{n \in L_1 \cap I_k} \lambda_n,\;\, \sum_{n \in I_k \setminus L_1} \lambda_n\; \text{ and }\; \sum_{n\in I_k}\lambda_n
$$
are all non-zero and have the same signs.
Moreover, $k$ must be an odd integer by part (4) of Proposition \ref{prop3.1}.

(iv)  By the previous part, the scalars:
$$
\e = \sum_{m < k} \sum_{n \in I_m} \lambda_n, \;\, \gamma = \sum_{n \in L_1 \cap I_k} \lambda_n\, \text{ and }\, \delta = \sum_{n \in I_k \setminus L_1} \lambda_n
$$
are all non-zero and have the same sign.

(v)  Since $L_1\in \mathcal{L}$, (ii) and part (2) of Proposition \ref{prop3.1} yield
\[\sum_{m < k} \bigl | \sum_{n \in I_m} \lambda_n \bigr |^2 + \gamma^2 \leq  \|x|L_1\|_J^2 = \bigl |\sum_{n \in L_1} \lambda_n \bigr |^2
= ( \gamma + \varepsilon )^2.\]

Setting $\rho = \sum_{m < k}  | \sum_{j \in I_m} \lambda_j |^2$, (iv) implies that $\varepsilon \delta > 0$ and $\rho + \gamma^2 \leq (\gamma + \varepsilon)^2$. Then, by Lemma \ref{l1}, $\rho + (\gamma + \delta )^2 < (\gamma + \delta + \e)^2$.
Hence
\[\sum_{m =1}^k \bigl| \sum_{n \in I_m} \lambda_n \bigr |^2 < \bigl |\sum_{n \in \cup_{m=1}^k I_m} \lambda_n \bigr|^2.\]

Thus $\{I_i\}_{i\in F}$ is not an $x$-norming partition, because $\mc Q=\{\cup_{m \leq k}I_m\}\cup \{I_m\}_{m >k}$ gives an estimate on the norm of $x$ strictly greater than that of $\{I_i\}_{i \in F}$. This is a contradiction that completes the proof.
\end{proof}

The $x$-norming partitions are not necessarily unique, but they have some regularity properties that we describe in the following corollaries.

\begin{corollary}\label{cor3.3}
Let $x\in J$ and let $\mc{I} = \{I_i\}_{i \in F}$ and $\mc{L}= \{L_j\}_{j\in G}$ be $x$-norming partitions.
Then, given $i\in F$ and $j \in G$, we have $I_i\subset L_j$ or $L_j\subset I_i$ or $I_i\cap L_j = \emptyset$.
\end{corollary}
\begin{proof}
The sets
$$
F_1 = \bigl\{ i\in F:\, \exists j\in G,\, \min I_i =\min L_j \bigr\}
$$
$$G_1 = \bigl\{ j\in G:\, \exists i\in F, \, \min I_i = \min L_j \bigr\}
$$
are non-empty since both contain $1$.
Let $i_0\in F_1$. Then there exists a unique $j_0\in G_1$ such that $\min I_{i_0}= \min L_{j_0} $. Hence, by Proposition \ref{p1},
$$
I_{i_0} \subset L_{j_0} \quad \text{or} \quad I_{i_0}\cap \supp(x) = \cup_{j\in Q} [L_j \cap \supp(x)]
$$
for some subinterval $Q$ of $G$.

In the first alternative $i_0$ and every $j\in G$ fulfill the conclusion, while in the second one for every $j\in Q$ we have $L_j \subset I_{i_0}$ and any other $L_k$ is disjoint from $I_{i_0}$ (Proposition \ref{p1}). The same holds for $j_0 \in G_1$.

Assume that $i_0 \in F\setminus F_1$ and let $i_1= \max\{i\in F_1: i<i_0\}$ and $j_1 \in G_1$ such that $\min I_{i_1} = \min L_{j_1}$. Then Proposition \ref{p1} yields
$$
L_{j_1}\cap \supp(x) = \cup_{i\in Q} [I_i \cap \supp(x)]\quad \text{and} \quad i_0 \in Q.
$$
Therefore $I_{i_0}$ is disjoint to each $L_j$ with $j\ne j_0$, and the result is proved. Similar arguments work in the case $j_0\in G\setminus G_1$.
\end{proof}

Let us see that there exists a joint refinement of any two $x$-norming partitions.

\begin{corollary}\label{cor3.4}
Let $x\in J$ and let $\mathcal{I}= \{I_i\}_{i \in F}$ and $\mathcal{L}= \{L_j \}_{j\in G}$ be $x$-norming partitions. Then there exists an $x$-norming partition $\mathcal{N} =\{N_p\}_{p\in H}$ such that

(1) $\mc{N} \subset \mc{I} \cup \mc{L}$.

(2) For every $i \in F$ and $j\in G$ there exist subintervals $H_i$ and $H_j $ of $H$ such that
$$
I_i\cap \supp(x)= \cup_{p\in H_i}  [N_p\cap \supp(x)] \text{ and }
$$
$$L_j\cap \supp(x)  =\cup_{p \in H_j} [N_p\cap \supp(x)].
$$
\end{corollary}
\begin{proof}
We define
$$
\mathcal{N}= \bigl\{ I_i : \exists j\in G, \, I_i \subset L_j \bigr\} \cup \bigl\{L_j: \exists i\in F,\, L_j \subset I_i \bigr\}
$$
and write $\mathcal{N}=  \{N_p\}_{p\in H}$. Clearly $\mc{N}$ satisfies (2). Let us see that $\mc{N}$ is an $x$-norming partition.

For $i\in F$, either $I_i\in \mc N$ or $I_i \cap \supp(x)  = \cup_{p\in H_i}[N_p  \cap \supp(x) ]$.
In the second case, for each $p \in H_i$ there exists $j_p \in G$ such that
$N_p  = L_{j_p}$. Since both $\mc I$ and $\mc L$ are $x$-norming, we have
$$
\bigl| \sum_{n\in I_i}x(n)\bigr|^2 = \sum_{p\in H_i} \bigl| \sum_ {n\in N_p}x(n) \bigr|^2.
$$
Hence $\sum_{p\in H}\bigl|\sum_{n\in N_p} x(n)\bigr|^2  = \sum_{i\in F} \bigl| \sum_{n\in I_i}x(n)\bigr|^2$, proving that $\mc N$ is $x$-norming.
\end{proof}

In Proposition \ref{prop3.5.1} we will show the existence of a finest $x$-norming partition for each non-zero $x\in J$.

 \begin{definition}
A vector $x\in J $ is said to be \emph{separated by partitions} if, given $n_1 < n_2$ consecutive elements of $\supp(x)$, there exists an $x$-norming partition $\mathcal{I}= \{I_i\}_{i \in F}$ such that $n_1$ and $n_2$ belong to different intervals of $\mathcal{I}$.
\end{definition}

The next result is a consequence of Corollary \ref{cor3.4}, and it is the basic ingredient in the proof of the characterization of the extreme points of $B_J$.

\begin{proposition}\label{prop2.1}
If $0\neq x\in J$ is separated by partitions then $\{\{n\} \}_{n\in \supp(x)}$ is an $x$-norming partition. Hence $\|x\|_J = \|x\|_2$.
\end{proposition}
\begin{proof}
Let $\supp(x)= (m_k)_{k\in F}$. We show by induction that for every $k\in F$ there exists an $x$-norming partition $\mathcal{I}_k$ that includes the intervals $\{m_1\}, \ldots, \{m_k\}$.


Any $x$-norming partition separating $m_1, m_2$ can be taken as $\mathcal{I}_1$. Assume that $k+1 < \sup F$ and $\mathcal{I}_k$ exists. Then we choose an $x$-norming partition $\mathcal{L}_{k+1}$ which separates $m_{k+1}, m_{k+2}$. By Corollary \ref{cor3.4} there exists a joint refinement $\mc I_{k+1}$ of $\mathcal{L}_{k+1}$ and $\mathcal{I}_k$ that satisfies the inductive conclusion. If $k+1= \sup F$ then $\mathcal{I}_k$ is an $x$-norming partition consisting of singletons. This completes the inductive proof.

If $F$ is finite, the result has been proved. If $F= \N$, we set $\mathcal{I}$ as the pointwise limit of $(\mathcal{I}_k)_k$. By Proposition \ref{P0}, $\mathcal{I}$ is the $x$-norming partition that we desire.
\end{proof}

\subsection{On the extreme points of $B_J$}
The main result of this subsection is the following characterization of the extreme points of $B_J$, denoted $Ext(B_J)$.

\begin{theorem} \label{T1}
Let $x \in J$. Then $x\in Ext(B_J)$ if and only if $\|x\|_J = \|x\|_2  =1$.
\end{theorem}

We postpone the proof of Theorem \ref{T1}.
First, we collect some facts about the extreme points of $B_J$.

\begin{proposition} \label{P02}
Every $x\in Ext(B_J)$ is separated by partitions, and for $n_1$ and $n_2$ consecutive members of $\supp(x)$, we have $x(n_1)x(n_2) <0$.
\end{proposition}
\begin{proof}
Suppose that $x= \sum_n a_n e_n \in S_J$ and $n_1 < n_2$ are consecutive elements of $\supp(x)$ not separated by partitions. Corollary \ref{C0} yields
$$
\bigl\| \sum_{n \leq n_1} a_n e_n \bigr \|_J^2 + \bigl\|\sum_{n \geq n_2} a_n e_n\bigr \|_J^2 < 1.
$$
Thus, by the continuity of the norm, there exists $t > 0$ small enough so that
\[ \bigl \| \sum_{n < n_1} a_n e_n  + (a_{n_1} + t)e_{n_1} \bigr \|_J^2 +
\bigl \|(a_{n_2} - t)e_{n_2} + \sum_{n > n_2} a_n e_n \bigr \|_J^2 < 1\, \text{ and } \]
\[ \bigl \| \sum_{n < n_1} a_n e_n  + (a_{n_1} - t)e_{n_1} \bigr \|_J^2 +
\bigl \|(a_{n_2} + t)e_{n_2} + \sum_{n > n_2} a_n e_n \bigr \|_J^2 < 1. \]

Now it is clear that
\[ \bigl \| \sum_{n < n_1} a_n e_n  + (a_{n_1} + t)e_{n_1} + (a_{n_2} - t)e_{n_2} + \sum_{n > n_2} a_n e_n \bigr \|_J^2 \leq  1\, \text{ and } \]
\[ \bigl \| \sum_{n < n_1} a_n e_n  + (a_{n_1} - t)e_{n_1} + (a_{n_2} + t)e_{n_2} + \sum_{n > n_2} a_n e_n \bigr \|_J^2 \leq  1. \]
Hence $x$ is the average of two different vectors in $B_J$, thus $x\notin Ext (B_J)$. 

Also, if $n_1$ and $n_2$ are consecutive members of $\supp(x)$, by the first part of the result there is an $x$-norming partition $\mc I$ and successive intervals $I_i, I_{i+1}$ in $\mc I$ such that $n_1 = \max I_i$ and $n_2 = \min I_{i+1}$. Applying parts (3) and (4) of Proposition \ref{prop3.1}, we get $x(n_1)x(n_2)< 0$.
\end{proof}

\begin{proof}[Proof of Theorem \ref{T1}]
Suppose that $\|x\|_J=\|x\|_2=1$ and $x = (y+z)/2$ for some $y,z\in B_J$. Since $\|u\|_2 \leq \|u\|_J$ for all $u \in J$, we get $y, z\in B_{\ell_2}$. Thus $x=y=z$ because $x\in Ext(B_{\ell_2})$. Hence $x\in Ext(B_J)$.

Conversely, for $x= \sum_{j=1}^\infty a_j e_j \in Ext(B_J)$ we have $\sum_{j=1}^\infty |a_j|^2 =1$.
Indeed, Proposition \ref{P02} yields that any $n_1< n_2$ in $\supp(x)$ are separated by an $x$-norming partition. Hence Proposition \ref{prop2.1} provides the result.
\end{proof}

\subsection{Further properties of $Ext(B_J)$}
In this subsection we point out a few direct corollaries of Theorem \ref{T1}. We start with the following result stated in \cite{B2}.

\begin{corollary} \label{C31}
The set $Ext(B_J)$ is norm-closed.
\end{corollary}
\begin{proof}
Since $\|x\|_2\leq \|x\|_J$, there is a continuous natural inclusion $i:J\to \ell_2$, and $Ext(B_J)= S_J\cap i^{-1}(S_{\ell_2})$.
\end{proof}

Given $x \in J$  with $\supp(x) = \{ m_1< \cdots < m_k < \cdots \}$ and a set $F= \{ n_1 < \cdots < n_k < \cdots\}\subset \N$, we denote by $x^F$ the element of $J$ satisfying $x(m_k) = x^{F} (n_k)$ for each $k \in \N$ and $suppx^F = F$.
Similarly, we define $x^F$ when $\supp(x)$ and $F$ are finite sets.

\begin{corollary}
For $x$ and $F$ as above, $x$ is in $Ext (B_J)$ if and only if so is $x^F$.
\begin{proof}
Since $(e_n)$ is spreading for $\|\cdot\|_J$, the map $e_{m_k} \to e_{n_k}$ induces an isometry between $ \overline {\langle(e_{m_k})_k \rangle}$ and $\overline {\langle(e_{n_k})_k \rangle}$. Hence $ \| x\|_J  = \| x^F \|_J$. And clearly $\| x\|_2  = \| x^F \|_2$.
\end{proof}
\end{corollary}

\begin{corollary} \label{C32}
Let $x\in S_J$.
\begin{enumerate}
\item  If $x\in Ext(B_J)$ then $\|x|_I\|_{J} = \|x|_I\|_2 $ for every interval $I$.
\item If $(I_i)_{i\in F}$ is an $x$-norming family and we set $b_i= \sum_{n \in I_i} x(n)$ for every $i \in F$,
then $\sum_{i\in F} b_{i} e_{k_i}\in Ext(B_{J})$ for every strictly increasing sequence $(k_i)_{i\in F}$ $\N$.
\end{enumerate}
\end{corollary}
\begin{proof}
(1) is a consequence of Proposition \ref{prop3.1}.

(2) We prove the case $k_i=i$ for all $i\in F$. The general case follows from Theorem \ref{T1} and the fact that
$(e_n)$ is a spreading basis for $J$.

Since $(I_i)_{i\in F}$ is an $x$-norming partition for $x\in S_J$, $\sum_{i \in F}|b_i|^2 =1$. Let $(E_i)_{i=1}^d$
be a finite sequence of bounded disjoint intervals of $\N$. Since $\cup_{i\in F}I_i\supset \supp (x)$,
$$
\sum_{r=1}^d \biggl | \sum_{i \in E_r} b_i \biggr |^2 = \sum_{r=1}^d \biggl | \sum_{i \in E_r} \sum_{j \in I_i} x(j)
\biggr |^2 \leq \|x\|_J^2.
$$
Then $y= \sum_{i\in F} b_i e_{k_i}\in J$ with $\|y\|_J \leq 1$, hence $\|y\|_J = \|y\|_2 =1$, and by Theorem \ref{T1}
we conclude that $y\in Ext(B_J)$.
\end{proof}

\begin{definition}
A vector $x \in J$ has \emph{non-positive remainder,} denoted NPR, provided that $\bigl | \sum_n x(n) \bigr |^2 \leq \sum_n |x(n)|^2$.

The vector $x$ has \emph{NPR hereditarily} if $x|_I$ has NPR for every interval $I$ of $\N$.
\end{definition}

Our next corollary is an immediate consequence of the first part of Corollary \ref{C32}.

\begin{corollary} \label{C33}
Let $x \in J$. Then $x\in Ext(B_J)$ if and only if $ x$ has NPR hereditarily and $\|x\|_2 =1$.
\end{corollary}

Next we give a characterization of the norm $\|\cdot \|_*$ that does not involve $\|\cdot\|_J$.

\begin{corollary} \label{C34}
Let $x^* \in J^*$. Then,
\[\|x^*\|_* = \sup \{ x^*(x): \, x \text{ has NPR hereditarily and } \|x\|_2 =1\}.\]
\end{corollary}

We close this subsection with a characterization of the extreme points of $B_{J_s}$, where $J_s$ is the James space
defined with the square variational norm in Formula (\ref{s-norm}).

\begin{corollary}\label{cor3.2}
Let $x= \sum_n a_n e_n\in J_s$ with $\|x\|_s= 1$. Then $x\in Ext(B_{J_s})$ if and only if the sequence $(a_n- a_{n+1})_n$
is $NPR$ hereditarily.
\begin{proof}
As we have mentioned, the map $ (c_n)_n \to (c_{n} - c_{n+1})_n$ is an onto isometry from $J_s$ to $J$.
Hence $(c_n)_n$ is an extreme point of $B_{J_s}$  if and only if $(c_{n} - c_{n+1})_n$ is an extreme point of $B_J$, and the result follows from Corollary ~\ref{C33}.
\end{proof}
\end{corollary}

\subsection{The finest $x$-norming partition}
Here we show that for every $0\neq x\in J$ there exists a finest $x$-norming partition which corresponds to the
maximal cost path in  \cite {B2}.

\begin{definition}
Let $x\in J$ and let $\mathcal{I} = \{ I_i\}_{i\in F}$ and $\mathcal{L}= \{L_j\}_{j\in G}$ be $x$-norming partitions.
We say that \emph{$\mathcal{I}$ refines $\mathcal{L}$} if for every $j\in G$ there exists a subinterval $F_j$ of $F$
such that $L_j \cap\supp(x)= \cup_{i\in F_j} [I_i\cap \supp(x)]$.
\smallskip

An $x$-norming partition $\mathcal{I}$ is said to be the \emph{finest $x$-norming partition} if it refines any
other $x$-norming partition.
\end{definition}

If $\mathcal{I} = \{I_i\}_{i\in F}$ and $\mathcal{L}= \{L_j\}_{j\in G}$ are $x$-norming partitions and $\mathcal{I}$
refines $\mathcal{L}$ then for every $i \in F$ and $j\in G $ either $I_i \subset L_j$ or $I_i \cap L_j = \emptyset$.

\begin{proposition} \label{prop3.5.1}
Every non-zero $x\in J$ admits a finest $x$-norming partition.
\end{proposition}
\begin{proof}
Let $x \in J$, $x\ne 0$. For $n\in \supp(x)$ we denote by $\mc{I}_n$ the family of all intervals $I$ such that
$n\in I$ and there exists an $x$-norming partition $\mc{I}$ with $I\in \mc{I}$.

Corollary \ref{cor3.3} yields that $(\mc{I}_n, \subset)$ is totally ordered; e.g. for $I_1, I_2 \in \mc{I}_n$
either $ I_1 = I_2$ or $ I_1 \subset  I_2$ or  $ I_2 \subset I_1$. Therefore for every $n\in \supp(x)$ there
exists a unique $I_n \in \mc{I}_n$ such that  $I_n \subset I$ for every $I\in \mc{I}_n$.

Observe that Corollary \ref{cor3.3} also implies that for $n \ne m \in \supp(x) $ either $I_n = I_m$ or
$I_n \cap I_m = \emptyset$. Therefore the family $\{I_n\}_{n\in \supp(x)}$ defines a partition of the support of $x$ which refines any $x$-norming partition. We reorder this family as $\mc{L} =\{ L_i\}_{i\in F} $ with $F$ an initial segment of $\N$ or $F=\N$ and for $i<j \in F$, $\max L_i < \min L_j$.

Using induction we show that $\mc{L}$ is an $x$-norming partition. More precisely for $k\in F$ we prove that $\{L_i\}_{i\leq k}$ is a  $x_k$-norming partition, where $ x_k = x|_{\cup_{i\leq k} L_i}$.

For $k = 1$  the result follows from Proposition \ref{prop3.1} (2). Assume the inductive assumption
proved for $k\in F$, and let $\mc{M} =\{ M_p\}_{p\in G}$ be an $x$-norming partition with
$L_{k+1} = M_{p_0}$. Then $ \{ M_p\}_{p < p_{0}}$ is an $x_k$-norming partition and by the inductive
assumption $\{L_i\}_{i\leq k}$ shares this property. Thus $\{L_i\}_{i\leq k+1} \cup \{M_p\}_{ p_0< p}$ is an $x$-norming partition and Proposition \ref{prop3.1} (2) yields that  $\{L_i\}_{i\leq k+1}$ is an $x_{k+1}$-norming partition, and the proof is complete. 
\end{proof}

\begin{remark}
The existence of the finest $x$-norming partition and Proposition \ref{P02} easily imply that for every
$x\in Ext(B_J)$ we have $\| x\|_J = \| x\|_2$.
Indeed,  $x$ is an extreme point if and only if its finest $x$-norming partition consists of singletons.
\end{remark}

\section{The norm of $J^{**}$ and the extreme points of $B_{J^{**}}$}\label{sect:J**}

Here we give an expression for the norm of $J^{**}$ and present a characterisation of the extreme points of the unit ball of $J^{**}$.

\subsection{A norming set for $J^{**}$}
In the sequel, $I_{\infty }^*$ and $I_{[k, \infty)}^*$ denote the functionals in $J^*$ that correspond to the set $\N$ and to the final segment of $\N$ starting from $k\in \N$.

The next result will be necessary later. It is well-known (\cite {A2}, \cite{CLL}), and it is also a consequence of the main result of \cite{AI}.

\begin{lemma}\label{lem4.2} Let $(x_n)_n$ be a normalized  block sequence in $J$ such that  $I_{\infty}^* (x_n) \to 0$. Then $(x_n)$ has a subsequence equivalent to the $\ell_2$ basis.
\end{lemma}

The following result is also well-known. We prove it for completeness. Here $(e_n^*)_n$ is the sequence of biorthogonal functionals of the basis of $J$.

\begin{lemma}\label{lem4.3}
\begin{enumerate}
\item The set $\{e_{n}^* : n\in \N\} \cup \{ I_{\infty}^* \}$ generates a norm-dense subspace of $J^*$.
\item $J^{*} = J_{*} \oplus  \langle I_{\infty}^{*}\rangle$, where $J_{*}$ is the closed subspace of $J^*$ generated by $(e_{n}^{*})$, which is a predual of $J$.
\item The canonical copy of $J$ in $J^{**}$ has codimension one.
\end{enumerate}
\end{lemma}
\begin{proof}
(1) If not, there exists $x^{**}\in S_{J^{**}}$ such that $\{e_{n}^* : n\in \N\} \cup \{ I_{\infty}^* \} \subset \ker  x^{**}$. Since $J^*$ is separable, by Goldstine's theorem there exists a normalized sequence $(x_k)_k$ such that $x_{k} \xrightarrow{w^*} x^{**}$.
Observe that for each $n\in \N$,  $x_{k}(n) \to 0 $. Hence, we may assume that $(x_k)_k $ is a block sequence and $I_{\infty}^* (x_k) \to 0$. So Lemma ~\ref{lem4.2} shows that $(x_k)_k$ has a subsequence equivalent to the $\ell_2$ basis, a contradiction.

(2) Since the basis $(e_n)_n$ is boundedly complete, the sequence $(e_{n}^*)_n$ generates the canonical embedding in $J^{*} $ of the predual $J_{*}$  of $J$. Clearly, $I_\infty^* \notin J_*$ and the first part of the proof yields the result.

(3) It is clear that $J^{**} = (J_{*} \oplus \langle I_{\infty}^{*}\rangle)^* \simeq J \oplus \mathbb{R}$. Since $J_*$ is the canonical copy of the predual of $J$ in $J^*$, in the above equality $J$ is the canonical copy of $J$ in $J^{**}$.
\end{proof}

\begin{definition}
We consider the set $\mathcal{D}\subset J$ given by 
\begin{equation}\label{eq:D-set}
\mathcal{D}= \left\{\sum_{i=1}^k a_i I_i^*: \{I_i\}_{i=1}^k \text{ disjoint finite intervals of }  \N, \, \bigl (\sum_{i=1}^k a_i^2 \bigr)^{1/2} \leq 1\right\},
\end{equation}
and we denote $\mathcal{D}_1= \overline{\mathcal{D}}^{w^*}$.
\end{definition}

\begin{remark}
As we said in Section \ref{sect:intro}, given a family $\{I_i \}_{i\in F}$ of intervals of $\N$, $F$ is a finite or infinite initial segment of $\N$, each $I_i$ is non-empty and $\max I_i < \min I_j$ for $i < j \in F$.  
It easily follows from the definition of $\|\cdot\|_J$ that $x^* = \sum_{i\in F} a_i I_i^*\in J^*$ with $\sum_{i=1}^k a_i^2 \leq 1$ satisfies $\|x^*\|_* \leq 1$. 
\end{remark}

The proof of the following result is easy.

\begin{proposition}
The set $\mathcal{D}$ is contained in $B_{J^*}$ and it is norming; i.e., for every $x\in J$, $\|x\|_J=
\sup\{|x^*(x)|: x^*\in \mathcal{D}\}$.
\end{proposition}


Next we state a special case of \cite[Theorem 13B, p. 74]{H:75}:

\begin{theorem}\label{thm:H}
Let $X$ be a real Banach space.
For a subset $A\subset B_{X^*}$ such that $A=-A$, the following conditions are equivalent:
\begin{enumerate}
\item The convex hull of $A$ is $w^*$-dense in $B_{X^*}$.
\item $A$ is norming.
\item $Ext(B_{X^*})\subset \overline{A}^{w^*}$.
\end{enumerate}
\end{theorem}




\begin{corollary}
The convex hull of $\mathcal{D}$ is $w^*$-dense in $B_{J^*}$ and $Ext(B_{X^*})$ is contained in  $\mathcal{D}_1$.
\end{corollary}


The following fact is a consequence of a Theorem of Bessaga and Pelczynski (\cite{BP}). We prove it for the sake of completeness.

\begin{lemma}\label{lem4.5}
The set $\mathcal{D}_1 \subset J^*$ norms $J^{**}$.
\end{lemma}
\begin{proof} We denote by $\| \cdot \|_{**} $ the norm of $J^{**}$ and we will show that

\[\|x^{**}\|_{**} = \sup _{x^* \in \mathcal{D}_1} |x^{**}(x^*)|\, \text{ for all }\, x^{**} \in J^{**}.\]

Let $x^* \in J^*$. Since $Ext(B_{X^*})\subset \mc{D}_1$, Choquet's theorem (\cite[Proposition 2.20]{FLP}) gives a probability measure $\mu$ on the $w^*$-Borel subsets of $B_{J^*}$ supported on $\mathcal{D}_1$ such that
$$
x^*(x)= \int_{\mathcal{D}_1} y^*(x) d\mu(y^*) \, \textrm{ for all }\, x \in J.
$$

Let $x^{**} \in J^{**}$. Since $J^*$ is separable, Goldstine's theorem gives a sequence $(x_n)$ in $X$ such that $\|x_n \| _J\leq \|x^{**}\|_{**} $ for all $n\in\N$ and $x^{**}= w^{*}$-$\lim_n x_n$.
By the dominated convergence theorem
$$
x^{**}(x^*) = \int_{\mathcal{D}_1} x^{**}(y^*) d \mu(y^*),
$$
which implies $|x^{**}(x^*)|\leq \sup\bigl\{|x^{**}(y^*)| : y^* \in \mathcal{D}_1\bigr\}$, proving the lemma.
\end{proof}

We denote by $\omega+1$ the set $\N\cup \{ \omega \}$ with the order that extends the usual one in $\N$.

\begin{remark}
(1) Since $I_{\infty}^{*}(e_n) =1$ for each $n \in \N$, $(e_n)$ is a non-weakly convergent, weakly Cauchy sequence in $J$. Therefore, there exists $e_{\omega} \in J^{**}$ such that $e_n \xrightarrow {w^*} e_{\omega}$.
\smallskip

(2) For $k\in \N$, $I_{ [k, \infty)}^{*} (e_{\omega}) =1$. Thus, $I_{[k,\infty )}^{*} (x^{**}) = \sum _{k\geq n} x^{**}(n) + x^{**}(\omega)$ for all $x^{**} \in J^{**}$.
\end{remark}

\subsection{The norm of $J^{**}$}

Here we consider $\mathcal{D}_1 $ as a subset of $J^{***}$, $w^*$ is the $w^*$-topology on $J^{***}$, and for an interval $I$ of $\omega + 1$,  $I^*\in J^{***}$ is the associated functional.

\begin{lemma}\label{lem4.6}
\begin{enumerate}
\item For every interval $I$ of $\omega+1$ there exists a sequence of intervals $(I_n)_n \subset \mathcal{D}_1$ such that $ I_n \to I$ pointwise. Hence $I_{n}^* \xrightarrow{w^{*}} I^*$.
\smallskip

\item $\overline{\mathcal{D}_{1}}^{w^{*}} = \{ \sum_{i\in F} \alpha_{i} I_{i}^{*} :  \sum_{i\in F} \alpha_{i} ^{2} \leq 1,\, (I_{i} )_{i \in F}\, \text{ disjoint intervals of }\, \omega+1\},$ where F is either a finite initial interval of $\N$ or $\N$ or $\omega+1$.
\end{enumerate}
\end{lemma}
\begin{proof}
(1) As we have said before, $I_{[k,\infty)}^*$ acting on $J^{**}$ coincides with $I^*$ where $I = [k, \infty) \cup \{ \omega + 1\}$. Therefore, the only intervals $I$ of $\omega + 1$ such that $I^{*} \notin  \mathcal{D}_1$  are $I_{0} = \{ \omega +1 \}$ and $I_{k} = [ k, \infty)$. For the first case we set $I_n =  [n, \infty) \cup \{\omega + 1\}$, $n\in \N$ and it easily follows that  $I_{n }^{*}  \xrightarrow{w^*} I_0^*$. In the second case we set  $I_n = [k,  n] $ , $ k < n $ and we have that $I_{n }^{*}  \xrightarrow{w^*} I_k^*$.
\smallskip

(2) Consider the functional $x^{***} = \sum _{ p\in \omega + 1} \alpha_p I_{p}^*$, where $ \sum _{ p\in \omega + 1} \alpha_{p}^2 \leq 1$ and $ \{I_p \}_p$ are non-empty disjoint finite intervals.
Clearly, $I_{\omega} = \{ \omega \}$. For $n \in \N$ we set
\begin{equation*}
x_{n}^{***} = \sum_{p =1}^{n} \alpha_p I_{p}^*+ \alpha_\omega I_{[k, \infty)}^{*}.
\end{equation*}
It is easy to check that $x_{n}^{***} \in \mathcal{D}_1$ and $x_{n}^{***} \xrightarrow{w^*} x^{***}$.

The remaining cases can be proved in a similar manner.
\end{proof}

For $x^{**} \in J^{**}$ and a disjoint family of intervals $\mathcal{I} = \{ I_i \}_{i \in F}$ of $\omega+ 1$ we denote
\begin{equation*}
\| x^{**} \|_{\mathcal{ I}} = (\sum_{i \in F} (\sum_{ p\in I_{i}} x^{**}(p))^{2})^ {1/2}.
\end{equation*}

The proof of the following result is easy.

\begin{lemma}\label{lem4.1} Let $\mathcal{I} = \{ I_i \}_{i \in F}$ be a disjoint family of intervals of $\omega+ 1$. Then
\begin{equation*}
\|x^{**}\|_{\mathcal {I}} = \sup\{ \sum _{i\in F} \alpha_{i} I_{i}^{*} (x^{**}): \sum_{i\in F} \alpha_{i}^{2} \leq 1 \} \text{ for each } x^{**} \in J^{**}.
\end{equation*}
\end{lemma}

Since $\mathcal{D}_1 \subset  \overline{\mathcal{D}_1}^{w*} \subset B_{J^{***}}$ and $\mathcal{D}_1$ norms $J^{**}$ (Lemma ~\ref{lem4.5}),  Lemmas~\ref{lem4.6} and ~\ref{lem4.1}   yield the next description of the norm of $J^{**}$:

\begin{proposition} \label{prop4.1}
For $x^{**}\in J^{**}$
\[
\| x^{**} \|_{**} = \sup\bigl\{ \bigl( \sum_{i\in F } |\sum_{n\in I_{i}} x^{**}(n)|^{2} \bigr)^{1/2} : \{I_{i} \}_{i\in F} \text{ disjoint intervals of }  \omega+1 \bigr \},
\]
where  $F$ is either a finite initial interval of $\N$, or $\N$, or $\omega + 1$.
\end{proposition}

\begin{remark} For $0\neq x^{**} \in J^{**}$, the $x^{**}$-norming partitions are defined in a similar manner as $x$-norming ones (Definition \ref{x-norm-part}). Similarly, for every $x^{**} \in J^{**}$ there exists an $x^{**}$-norming partition satisfying the properties of Propositions \ref{prop3.1} and \ref{p1}.
\end{remark}

\begin{proposition}
The family $(e_p)_{p\in \omega + 1} $ is a bimonotone transfinite basis for $J^{**}$.
\begin{proof}
For $ x= \sum_{i=1}^{n} \alpha_{i} e_{i} + \alpha_{\omega} e_{\omega} $ and $y=  \sum_{i=1}^{n} \alpha_{i} e_{i}$ we have $\|y\|_{**} \leq \|x\|_{**}$.

In fact, since $y\in \langle(e_{n})_{\in\N} \rangle$, its norm is determined by a finite sequence of disjoint finite intervals of $\N$. This family gives the same evaluation on $x$. The remaining cases are treated in a similar manner. We conclude that the basis $(e_p)_{p\in \omega +1}$ is bimonotone for $\langle(e_{p})_{p \in \omega+1}\rangle$, and this property is extended to $J^{**}$.
\end{proof}
\end{proposition}

\begin{remark} In the case of $J$ we first  defined the norm and then we found the norming set $\mathcal{D} \subset J^{*}$. For $J^{**}$ we followed a reverse way. First we proved  that the set $\mathcal{D}_1$ norms $J^{**}$, and then we used $\overline{\mathcal{D}_1}^{w^*}\subset J^{***}$ to define the norm of $J^{**}$.
\end{remark}

The dual spaces of quasireflexive spaces were studied by  Bellenot in \cite{B}.

\begin{remark}
(1) In the fourth dual of $J$, denoted by   $J^{2(**)}$, the basis $(e_n)_n$ of $J$ is not  $w^*$-converging to $e_\omega$ since $I_{\{ \omega\}}^*\in J^{***}$ separates $(e_n)_n$ from $e_\omega$, but $(e_n)_n$ remains non-trivially weak Cauchy. Hence $J^{2(**)}$ is generated by $(e_n)_{n\in\N}$, $e_\omega$ and $e_{\omega + 1}$ in the natural order, where $e_{\omega + 1}$ is the image of $e_\omega \in J^{**}$ by the canonical embedding into $J^{2(**)}$, and the $w^*$-limit of $(e_n)_n$ is the new $e_\omega$. Therefore $\{e_n: n\in \N\}\cup \{ e_{\omega+ 1}\} $ generates the canonical copy of $J^{**}$ in  $J^{2(**)}$.

The norm of $J^{2(**)}$ is defined using the intervals of $\omega + 2$ as we did for  $J^{**}$. This argument can repeated for all the spaces $J^{n(**)}$.

(2) The higher even duals of $J$ are more easily described if we view them as subspaces of $J\Q$,  where $\Q$ is the set of rational numbers. Namely, $J\Q$ is the completion of $c_{00}(\Q)$ endowed with a James type norm defined as follows:

For $x = (x(p))_{p\in \Q} \in c_{00}(\Q)$ we set
\begin{equation*}
\| x \|_{J\Q} = \sup \bigl \{\bigl ( \sum_{i=1}^{n} |\sum _{p \in I_{i} }x(p) |^2 \bigr) ^{1/2}\, :\, \{ I_i \}_{i } \text{ disjoint intervals of } \Q \bigr \}.
\end{equation*}

For every $n \in \N$, the space $J^{n(**)}$ is isometric to any subspace of $J\Q$ generated by $( e_{t_k})_{k\in \N} \cup (e_{s_m})_{m= 1}^n$ with $(t_k)_k \subset \Q$ strictly increasing and $(s_m)_{m=1}^n \subset \Q$  strictly decreasing with  $\lim_{k} t_k < s_n$. Moreover the subspace generated by $( e_{t_k})_{k\in \N} \cup ( e_{s_m})_{m= 1}^l , 1\leq l < n$ is isometric to the natural embedding of $J^{l(**)}$ into  $J^{n(**)}$.

The space $J^{\omega(**)}$ is the limit of $J^{n(**)}$. It is isometric to any subspace of $J\Q$ generated by two sequences $(e_{t_k})_{k\in \N}, (e_{s_m})_{m\in \N}$ with $(t_k)_k \subset \Q$  strictly increasing, $(s_m)_{m\in \N} \subset \Q$ strictly decreasing and  $\lim_{k} t_k < \lim_{m} s_m$.

These arguments can be extended inductively to $J^{\alpha (**)}$ for $ \alpha \leq \omega^2$.
\end {remark}

\subsection{The extreme points of $B_{J^{**}}$}
The proof of the following result is similar to the one we did for $Ext(B_J)$. So  it is left to the interested reader.

\begin{theorem}\label{th3.1} Let $x^{**} = \sum_{j=1}^\infty a_j e_j + \ a_{\omega} e_{\omega}\in S_{J^{**}}$. Then $x^{**}\in Ext(B_{J^{**}})$ if and only if $\|x^{**} \|_2=1$.
\end{theorem}

The  properties of the extreme points of $B_J$ described in Subsection 3.3, remain valid for the extreme points of $B_{J^{**}}$. In particular the following analogue of Corollary \ref{C32} is proved with similar arguments.

\begin{corollary} \label{cor4.1}
Let $x^{**}\in S_{J^{**}}$.
\begin{enumerate}
\item  If  $x^{**}\in Ext(B_{J^{**}})$ then $\|x^{**}|_I\|_{**} = \|x^{**}|_I\|_2 $ for every interval $I$.
\item If $(I_i)_{i\in F}$ is an $x^{**}$-norming family and  $b_i = \sum_{n\in I_i} x^{**}(n)$ for each $i \in F$, then $\sum_{i\in F} b_{i} e_{k_i}\in Ext(B_{J^{**}})$ for every strictly increasing sequence $(k_i)_{i\in F}$ in $\N$.
\end{enumerate}
\end{corollary}

\section{The extreme points of $B_{J^*}$}\label{sect:J*}
Here we provide a complete description of the extreme points of $B_{J^*}$. In this section, $F$ is a finite initial segment of $\N$ or $\N$.

Since $\mathcal{D}_1$ is $w^*$-metrizable, standard arguments yield the next result.

\begin{lemma}\label{lem5.1} The set $\mathcal{D}_1$ is the union of the sets
$$
\biggl  \{\sum_{i \in F} \alpha_i I_i^* : \{ I_i\}_{i\in F} \text{ disjoint finite intervals and } ( \sum _{i \in F} | \alpha_i|^2)^{1/2} \leq1 \biggr \}\quad \text{and}
$$
$$
\biggl  \{\sum_{i=1}^n \alpha_{i }I_{i}^* +\alpha_{n+1}I_{[k, \infty)}^*: \{I_i\}_{i\in F} \text{ disjoint finite intervals, }
(\sum_{i=1} ^{n+1}| \alpha_i|^2)^{1/2} \leq1 \biggr\},
$$
where $\max I_n < k$ in the second set.
\end{lemma}

\begin{remark}
\quad (1) Since for $x^* =\sum _{i\in F} \alpha_i I_i^* \in \mathcal{D}_1$ we have  $\|x^*\|_* \leq (\sum_{i\in F}| \alpha_i|^2)^{1/2}$, when all the segments $I_i$ are finite $x^* =\sum_F \alpha_{i \in F}I_{i}^*\in J_*$, the predual of $J$.

(2) For every $x^* \in \mathcal{D}_1 \cap J_*$ there exists $x\in S_J$ such that $x^*(x)= 1$.
%
\end{remark}

Our next result connects the set $\mathcal{D}_1$ with $Ext(B_J)$.

\begin{proposition}\label{prop5.1}
Let $x^*= \sum_{i\in F} \alpha_iI_i^* \in\mathcal{D}_1$. Then $\|x^*\|_*= 1$ if and only if $\sum_{i\in F} \alpha_i e_i$ is an extreme point of $B_J$.
\end{proposition}
\begin{proof}
Suppose $\|x^*\|_*=1$ and take $x^{**} \in S_{J^{**}}$ such that $x^{**} (x^*) = 1$ and set $\beta_i = I_{i}^{*}(x^{**})$ for $i \in F$.
Since $\sum_{i\in F} \alpha_i \beta_i = 1$, $\sum_{i\in F}| \alpha_i|^2\leq 1$ and $\sum_ {i\in F} |\beta_i|^2 \leq 1$ we conclude $\sum_ {i\in F} |\beta_i|^2= 1$. So $\{I_i\}_{i\in F}$ is an $x^{**}$-norming family and $\alpha_i = \beta_i$ for all $i\in F$. Thus, Corollary \ref{cor4.1} (2) yields $\sum_{i\in F} \alpha_i e_i \in Ext(B_J)$.
\medskip

Conversely, assume that $x^*= \sum_{i\in F} \alpha_{i }I_{i}^* \in \mathcal{D}_1$ and $x= \sum_{i\in F }\alpha_{i }e_i\in Ext(B_J)$. We define $y=  \sum_{i\in F} \alpha_i e_{n_i}$ where $n_i = \min I_i$  for all $i\in F$. Then $\|y\|_J= (\sum_{i\in F} |\alpha_i|^2)^{1/2}=1$ and $x^*(y) = 1$. Hence $\|x^*\|_*=1$.
\end{proof}

Proposition \ref{prop5.1} shows that the extreme points of $B_{J^*}$ are among the vectors $x^* \in \mathcal{D}_1$ whose coefficients form an extreme point of $B_J$.
The next result shows that not all these $x^*\in \mathcal{D}_1$ are in $Ext(B_{J^*})$.

\begin{proposition} \label{prop5.2}
Let $x= \sum_{n\in F}\alpha_n e_n\in Ext(B_J)$ and let $\{I_n\}_{n \in F}$ be a (disjoint) family of intervals of $\N$.
If $\cup_{n\in \supp(x)}I_n$ is not an interval of $\N$, then $x^{*}= \sum_{n\in F} \alpha_{n} I_{n}^*$ is not an extreme point of $B_{J^*}$.
\begin{proof}
By our assumption there are successive points $k< m$ in $\supp(x)$ such that $\max I_k +1< \min I_m$. Setting $I_{k,m}= (\max I_{k}, \min I_{m})$, we have
$$
I_{k, m} \cap (\cup_{n\in \supp(x)} I_n) = \emptyset.
$$

We set  $J_{k} = I_{k} \cup I_{k, m}$ , $J_{m} = I_{m} \cup I_{k, m}$  which are intervals of $\N$. Notice that if there exists $n\in F$ such that $I_{k, m} \cap I_n \ne \emptyset$ then $\alpha_n = 0$. We define
\begin{equation*}
y^{*} = \sum_{n\in F,\\n \neq k}\alpha_{n} I_{n}^{*} + \alpha_{k}J_{k}^{*} \quad \textrm{and} \quad   z^{*} = \sum_{n\in F,\\n \neq m}\alpha_{n}I_{n}^{*} + \alpha_{m}J_{m}^{*}.
\end{equation*}

From our assumptions, $\alpha_{k}, \alpha_{m} \neq 0$. Since $x\in Ext(B_J)$ and $k, m$ are successive in $\supp(x)$, they must have alternate signs. Assuming
$\alpha_m >0$ we get $0 \in ( \alpha_k, \alpha_m)$. Hence $\lambda \alpha_k + (1-\lambda)\alpha_m = 0$ for some $0<\lambda < 1$.  It is easy to check that $x^* = \lambda y^{*} + (1- \lambda)z^*$ with $y^{*}\neq z^{*}$.
\end{proof}
\end{proposition}

Propositions \ref{prop5.1} and \ref{prop5.2} yield the following result.

\begin{corollary}\label{cor5.1} Let $x^*\in Ext(B_{J^*})$. Then $x^*= \sum_{i\in F} \alpha_i I_i^*$ with $\alpha_i\neq 0$ for each $i\in F$, $x= \sum_{i\in F} \alpha_i e_i\in Ext(B_J)$ and $\cup_{n\in \supp(x)}I_n$ is an interval of $\N$.
\end{corollary}

It is natural to ask if Corollary \ref{cor5.1} gives a characterization of $Ext(B_{J^*})$. Next, we give a positive answer to this question.

\begin{lemma}\label{lem5.2} Let $x^*, y^* \in \mathcal{D}_1$ satisfying the following conditions:
\begin{itemize}
\item[(a)] $x^* = \sum_{i\in F} \alpha_i I_{i}^*$ with $\alpha_i \ne 0$ for every $i\in F$ and $w= \sum_{i\in F} \alpha_i e_i\in Ext(B_J)$.
\smallskip

\item[(b)] $y^*= \sum_{j\in G}\beta_j L_j^*$ with $\beta_j \ne 0$ for every $j\in G$, and for every $z =\sum_{i\in F} \alpha_i e_{k_i} $ with $k_i\in I_i$ we have $y^{*}(z)= 1$.
\end{itemize}
\smallskip

Then the following assertions hold:
\begin{enumerate}
\item For $x= \sum_{i\in F}\alpha_i e_{n_i}$ with $n_i= \min I_i$ for $i\in F$, $\{L_j\}_{j\in G}$ is an $x$-norming family.
\item Setting $H_j= \{i\in F: n_i\in L_j \}$ and $k_j = \min H_j$ for $j\in G$, we have
\[ \beta_j = \sign(\alpha_{k_j} )\bigl(\sum_{i\in H_j}| \alpha_{i}|^2 \bigr)^{1/2}\, \text{ for every }j \in G.\]
\item $\cup_{i \in H_j} I_i \subset L_j $ for every $j\in G$.
\item If, additionally, $\cup_{i\in F}I_i$ is an interval of $\N$ then $\cup_{j\in G} L_j$ is also an interval and for every $j\in G$ with $1<j< \sup G$ we have $L_j= \cup_{i\in H_j} I_i$.

Moreover, if there exists $ m \in \supp(y^*)$  with $\max \supp(x^*) < m$, then we have
$\sign(y^*(e_m)) =\sign(\alpha_d)$, where $d = \max \supp(w)$.
\end{enumerate}
\end{lemma}
\begin{proof}
(1) If $x= \sum_{i\in F}\alpha_i e_{n_i}$ then our assumptions yield $x^*(x)= y^*(x)= 1$. Then
$$
y^*(x) = \sum_{j\in G} \beta_j \sum_{i\in H_j} \alpha_i  = 1, \, \sum_j |\beta_j|^2 \leq 1, \, \textrm{ and }\, \sum_j | \sum_{i\in H_j} \alpha_i |^2 \leq 1,
$$
which implies $\sum_j|\sum_{i\in H_j} \alpha_i|^2= 1$; hence  $\{L_j\}_j$ is an $x$-norming family.
\smallskip

(2)  From the previous part we conclude that $\beta_j = \sum_{i\in H_j} \alpha_i$ for all $ j\in G$.
Furthermore, since $x$ is an extreme point and $\{ L_j\}_{j\in G}$ is an $x$-norming partition,
\[ |\sum_{ i\in H_j } \alpha _i |^2 =  \sum_{i\in H_{j}} | \alpha_{i}|^2  \quad \text{and} \quad \sign(\sum_{i\in H_j} \alpha_i)= \sign(\alpha_{k_j}).\]
Therefore, $\beta_j= \sign(\alpha_{k_j}) \bigl(\sum_{i\in H_j}| \alpha_{i}|^2 \bigr)^{1/2}\, \textrm{ for every }j \in G.$
\smallskip

(3) Assume that there exists $j_0\in G$ and $i_0\in H_{j_0}$ such that $I_{i_0} \backslash L_{j_0} \ne \emptyset$.
Observe that the $i_0$ is necessarily unique.
To derive a contradiction  we choose $m_{i_0} \in I_{i_0} \backslash L_{j_0}$, we set $Q_{j_0} = H_{j_0} \backslash \{i_0\}$ and
$$
z = \sum_{i \notin H_{j_0}} \alpha_{i }e_{n_i} + \sum_ { i\in Q_{j_0}} \alpha_{i} e_{n_i}+ \ \alpha_{i_0}e_{m_{i_0}}.
$$

From our assumptions, it follows $y^*( z) = 1$, and arguing as in (2) with $z$ instead of $x$  we conclude that
\[ |\beta_{j_0} | =\bigl (\sum_{i\in Q_{j_0}}| \alpha_{i}|^2 \bigr)^{1/2}. \]
 On the other hand, from (2) we have
 \[ |\beta_{j_0}| =\bigl(\sum_{i\in H_{j_0}}| \alpha_{i}|^2 \bigr)^{1/2}\, \text{ and also }\, \sum_{i\in H_{j_0}}|\alpha_i|^2 \ne \sum_{i\in Q_{j_0}}|\alpha_i|^2,\]
which yields a contradiction, completing the proof of (3).

(4) The first part is a consequence of (3) and the fact that $\cup_{i\in \supp(w)}I_i$ is an interval.
For the second part, since $\{ L_j\}_{j\in G}$ is an $x$-norming partition, Proposition \ref{prop3.1} (3) yields that for every $j\in G$  and $p_j = maxL_j \cap \supp(x)$
$$
\sign\bigl(\sum_{i\in H_j} \alpha_i \bigr)= \sign (\alpha_{p_j}),
$$
and by (2), $\sign(\beta_j)= \sign(\alpha_{p_j}$).
To finish the argument, observe that $m\in L_{j_0}$ with $j_0 = \max G$. Hence $\sign(y^*(e_m))= \sign(\beta_{j _0})= \sign(\alpha_d)$.
\end{proof}

\begin{remark}\label{rem5.1}
Suppose that $x^*= \sum_{i\in F}\alpha_i I_i^* \in \mathcal{D}_1$,
$x = \sum_{i\in F} \alpha_i e_i\in Ext(B_J)$ and $\cup_{i\in \supp(x)} I_i$ is an interval of $\N$. Then $F= \supp(x)$.
Thus, in the sequel, we may assume that $\supp(x) = F$ with $F$ an initial segment of $\N$ or $\supp(x) =\N$.
\end{remark}

In our next lemma we adopt the notation used in the statement of Lemma \ref{lem5.2}.

\begin{lemma}\label{lem5.4} Let $x^* = \sum_{i\in F }\alpha_i I_i^* \in \mathcal{D}_1$ with $\alpha_i \neq 0$ for every $i \in F$. Suppose that
$$
w =\sum_{i\in F}\alpha_i e_i\in Ext(B_J)\, \text{ and }\, \cup_{i\in F} I_i\, \text{ is an interval of }\, \N.
$$

Then there exists a  $w^*$-closed face $E$ of $B_{J^*}$ such that $x^*\in E$, and every
$y^*=\sum_{j\in G} \beta_j L_j^* \in E\cap\mathcal{D}_1 $ with $\beta_j\ne 0$ for $j\in G$ satisfies the following conditions:
\begin{enumerate}
\item  $\{L_j\}_{j\in G}$ is an $x$-norming family for $x= \sum_{i\in F} \alpha_i e_{n_i}$ with $n_i = \min I_i$.
\item  For every $j \in G$, if $H_j = \{i : n_i \in L_j \}$  and $k_j = minH_j $ then
$$
\beta_j= \sign(\alpha_{k_j})\bigl( \sum_{i\in H_{j}} | \alpha_{i}|^2 \bigr)^{1/2}.
$$
\item $ \cup_{ j\in G} L_j$ is an interval of $\N$, $\cup_{i\in H_j} I_i\subset L_j$ for each $j\in G$  and $ \cup_{ i\in H_j} I_i  = L_j$ for each $1< j < \sup G$.
Moreover, if $\max \supp(x^*)< \max\supp(y^*)$  then $\sign(y^*(e_m))= \sign(\alpha_d)$, where $d= \max \supp(w)$.
\end{enumerate}
\end{lemma}
\begin{proof}
We set $S = \bigl\{\sum_{i\in F}\alpha_i e_{m_i}\, :\, m_i \in I_i \textrm{ for all } i\in F\bigr\}$, and for each $s\in S$ we define
 \[E_s = \bigl\{ x^* \in B_{J^*} : x^*( s) = 1 \bigr\}.\]
Then $E_s$ is a $w^*$-closed face of $B_{J^*}$ and $x^* \in E_s$. Hence $ E = \cap_{s\in S} E_s$
is a non-empty $w^*$-closed face and $x^*\in E$.

If $y^* \in E\cap \mathcal{D}_1$ and $s\in S$ then $x^*(s)= y^*(s) =1$. Thus, for all $y^*\in E\cap \mathcal{D}_1$ the pair $x^*, y^*$ satisfies the conditions of Lemma \ref{lem5.2}, implying the desired properties of $y^*$.
\end{proof}

\begin{proposition} \label{prop5.3}
Let $x^* = \sum_{i\in F}\alpha_i I_i^*\in \mathcal{D}_1$ with  $\alpha_i \neq 0$ for all $i \in F$. Suppose that $x = \sum_{i=1}^n \alpha_ie_i\in Ext(B_J)$ and $\cup_{i \in F} I_i$ is an interval of $\N$. Then $x^*\in Ext(B_{J^*})$.
\end{proposition}

\begin{proof}
 Let $E$ be the face of $B_{J^*}$ defined in Lemma \ref{lem5.4}.  Then $Ext(E) \subset Ext(B_{J^*})\subset \mathcal{D}_1$, hence $\overline{co}^{w^*} ( \mathcal{D}_1 \cap E) = E.$

Since $\mathcal{D}_1 \cap E$ is $w^*$-compact, the easy version of Choquet's theorem (\cite[Theorem 2.21]{FLP}) yields that for $z^* \in E $ there exists a (not necessarily unique) probability measure $\mu_{z^*}$ with $ \supp(\mu_{z^*} ) \subset   \mathcal{D}_1 \cap E $ which represents $z^*$. Namely, for every $w \in J$,
$$
z^*(w) = \int_{E} y^*(w) d \mu_{z^*}(y^*).
$$

We denote by $\mu_{x^*}$ a probability measure  with  $ \supp(\mu_{x^*} ) \subset   \mathcal{D}_1 \cap E $ representing $x^*$.
We have $\mu_{x^*} = \delta_{x^*}$, hence $x^* \in Ext(E)$. Otherwise $x^* = (z^* + w^*)/2$ with $z^* \ne w^*$ and if  $\mu_{z^*}$, $\mu_{w^*}$  represent $z^*, w^*$ then $\mu_{x^*} = (\mu_{z^*} + \mu_{w^*})/2 \ne \delta_{x^*}$.

Let $y^*= \sum_{j\in G} \beta_j L_j^* \in \mathcal{D}_1 \cap E$ with $\beta_j\ne 0$ for all $j \in G$. Part (3) of Lemma \ref {lem5.4} yields $\min L_1 \leq \min I_1$.
We set
$$
W_0= \biggl\{ y^*\in \mathcal{D}_1 \cap E:  y^* = \sum_{j\in G} \beta_j L_j^* \text{ and } \min L_1 < \min I_1 \biggr\}\, \text{ and}
$$
$$
V_0 = \biggl\{y^* \in \mathcal{D}_1 \cap E :  y^* = \sum_{j\in G} \beta_j L_j^* \, \text{ and }\, \min L_1= \min I_1\biggr\}.
$$

Since for $y^*= \sum_{j\in G} \beta_j L_j^* \in   \mathcal{D}_1 \cap E$ we have $\min L_1 \leq \min I_1$, we conclude
$$
W_0 \cup V_0 = \mathcal{D}_1 \cap E.
$$

Let us show that $\mu_{x^*}(W_0) = 0$. Indeed, assume that $ W_0 \ne \emptyset$ and set $ m= \min I_1 -1$.
Then for $y^*= \sum_{j\in G} \beta_j L_j^*\in W_0$, Lemma \ref{lem5.4} yields
$$
y^*(e_m) = \beta_1, \sign(\beta_1)= \sign(\alpha_1), \, |\beta_1| \geq |\alpha_1|.
$$
Assuming that $\mu_{x^*} (W_0) > 0$, since $y^*(e_m)= 0$ for $y^* \in V_0$, we have
 \begin{align*}
 0 < \biggl| \int_{W_0} y^*(e_m)d\mu_{x^*} \biggr|=\biggl| \int_{W_0} y^*(e_m)d\mu_{x^*}  + \int_{V_0} y^*(e_m)d\mu_{x^*} \biggr|=  \\
  \biggl| \int_{   \mathcal{D}_1 \cap E } y^*(e_m)d\mu_{x^*} \biggr| =  \biggl| \int_{E} y^*(e_m)d\mu_{x^*}\biggr| =| x^*(e_m)| = 0,
 \end{align*}
a contradiction. Hence $\mu_{x^*}(W_0) = 0$. Moreover, we have $\mu_{x^*}(V_0)= 1$ because $\mathcal{D}_1 \cap E \setminus W_0 = V_0$.

For $k\in F$ with $k< \sup F$ we set
\[ V_k = \biggl\{ y^* \in   \mathcal{D}_1 \cap E :  y^* = \sum_{ i\leq k} \alpha_i I_i^* + \sum_{ \substack {j\in G,\, k<j}} \beta_j L_j^* \biggr\},\]
where $G$  is an initial interval of $\N$ that depends on $y^*$, and we inductively show that $\mu_{x^*} (V_k) = 1$.
We first prove the inductive assumption for $k=1$ and $ 1< \sup F$.

Observe that $ V_1 \subset V_0 $.  If $ y^*=  \sum_{j\in G} \beta_j L_j^* \in V_0$ with $|H_1| = 1$  (see Lemma \ref{lem5.2} for the definition of $H_1$) then, since $1<\sup F$, we have $L_1= I_1$. Indeed, $I_2\subset L_2$, $I_1\cup I_2 $ is an interval and $\min L_1 = \min I_1$. This yields that $y^* \in V_1$,  hence
$$
W_1= V_0\setminus V_1= \bigl\{y^*\in V_0 : y^*=  \sum_{j\in G} \beta_j L_j^* \text{ and } |H_1|> 1 \bigr\}.
$$

We will show that $\mu_{x^*}(W_1) = 0$. We fix $m\in I_1$. Then $y^* = \sum_{j\in G} \beta_j L_j^*\in W_1$ satisfies
$$
y^*(e_m) = \beta_1,\, \sign(\beta_1)= \sign (\alpha_1)\, \text{ and } \, |\beta_1| \geq (|\alpha_1|^2 +|\alpha_2|^2)^{1/2}
$$
by Lemma \ref{lem5.4} (2).
Also for $y^* \in V_1$ we have that $y^*(e_m) = \alpha_1$. Therefore, assuming that $\mu_{x^*} (W_1) > 0$, we have
$$
\int_{W_1} y^*(e_m)d\mu_{x^*}= \gamma\mu_{x^*} (W_1),\, \sign(\gamma) = \sign(\alpha_1) \text{ and } |\gamma| > |\alpha_1|.
$$
Hence
\begin{eqnarray*}
|x^*(e_m)| &=&  \biggl| \int_{W_1} y^*(e_m) d\mu_{x^*}+ \int_{V_1} y^*(e_m)d\mu_{x^*} \biggr|\\
&=& \bigl|\gamma \mu_{x^*}(W_1)+ \alpha_1 \mu_{x^*} (V_1) \bigr| > |\alpha_1|,
\end{eqnarray*}
which gives a contradiction since $x^*(e_m) = \alpha_1$. Thus $\mu_{x^*}(V_1)= \mu_{x^*}(V_0) = 1$ and the proof for $k=1$ is complete.

Let $k\in F$, $1< k <\sup F$ and assume that $\mu_{x^*}(V_{k-1}) = 1$. To show that $\mu_{x^*}(V_{k}) = 1$ we follow the same steps as in the case $k=1$ with $V_{k-1}$ in the place of $V_0$. This completes the inductive proof.

If $x^* = \sum_{i=1}^\infty \alpha_iI_i^* $ then $\mu_{x^*} = \delta_{x^*}$ since $\cap_{k\in \N} V_k =\{x^*\} $ and for all $k\in  \N$, $ \mu_{x^*}(V_k) =1$.\\
For  $x^* = \sum_{i=1}^d \alpha_iI_i^* $ with $I_d$ an infinite subinterval of $\N$, then we have proved by induction that $\mu_{x^*}(V_{d-1}) = 1$ and
   \[ V_{d-1} =  \biggl\{ y^* \in   \mathcal{D}_1 \cap E :  y^* = \sum_{ i\leq d-1} \alpha_i I_i^* + \beta_d L_d \biggr\}. \]

It follows from part (3) of Lemma \ref{lem5.4} that $I_d \subset L_d$. Since $\cup_{i\in F} I_i$ is an interval, we have $I_d = L_d $ and $\beta_d = \alpha_d$. Hence $V_{d-1} = \{x^*\}$, and thus $\mu_{x^*} = \delta_{x^*}$.

It remains the case when $x^*= \sum_{i=1}^d \alpha_iI_i^*$ and $I_d$ is a finite interval. In this case $V_{d-1} =\{x^*\} \cup W_d$, where
$$
W_d = \bigg\{ y^* \in V_{d - 1}: y^* = \sum_{i=1}^{d-1} \alpha_i I_i^* + \alpha_d L_d \quad \max I_d< \max L_d  \biggr\}.
$$

Inductively we have shown that $\mu_{x^*}(V_{d-1}) = 1$ . Finally, we prove that $\mu_{x^*}(W_d) = 0$.
Indeed, set $m = \max I_d + 1$ and observe that  for $y^* \in W_d$ , $y^*(e_m) =\alpha_d \ne 0$, while $x^*(e_m) = 0$. Therefore,
$$
0= x^*(e_m) = \int_{V_{d-1}}y^*(e_m)d\mu_{x^*} =   \int_{W_{d}}y^*(e_m)d\mu_{x^*}= \alpha_d \mu_{x^*}(W_d),
$$
which yields $\mu_{x^*}(W_d) = 0$.
Hence $ \mu_{x^*}( \{x^*\}) = \mu_{x^*} (V_{d-1}) =1$, implying $\mu_{x^*} = \delta_{x^*}$.
\end{proof}

Corollary \ref{cor5.1} and Proposition \ref{prop5.3} imply the main characterization of $Ext(B_{J^*})$:

\begin{theorem}\label{thm5.1} 
A vector $x^*\in J^*$ is in $Ext(B_{J^*})$ if and only if $x^* = \sum_{i\in F} \alpha_{i} I_i^*$ with  $\alpha_i\neq 0$ for every $i\in F$ and the following properties are satisfied:
\begin{enumerate}
\item $x= \sum_{i\in F}\alpha_i e_i\in Ext(B_J)$ and 
\item $\cup_{i\in F}I_i$ is an interval of $\N$.
\end{enumerate}
\end{theorem}

In contrast to $Ext(B_J)$ (see Corollary \ref{C31}), $Ext(B_{J^*})$ is not norm-closed as the next example shows.

\begin{example}\label{ex5.1}
We set
 \[ x^* = \frac{1}{\sqrt{2}}( e_1^* - e_4^*)\, \text{ and }\, x_m^* =\frac{m}{\sqrt {2(m^2 + 1)}}( e_1^* - \frac{1}{m} e_2^* + \frac{1}{m} e_3^* -e_4^*)\, \text{ for }\, m\in \N. \]

It is easy to check that $x^*, \, (x_m^*)_m$ are norm-one vectors such that
    \[ x^* \notin Ext(B_{J^*}),\, (x_m^*)_m \subset Ext(B_{J^*})\, \text{ and }\, \|x_m^*-x^*\|_*\to 0.\]
\end{example}

This example raises the question whether $\overline{Ext(B_{J^*})}^ {\| \cdot \|_* }= \mathcal{D}_1 \cap S_{J^*}$. We show that this is not the case.

\begin{definition}\label{def5.1} Let $x^* \in \mathcal{D}_1$ with $x^*= \sum_{i\in F}\alpha_i I_i^*$, where $F$ is an initial interval of $\N$ and $\alpha_i \ne 0$ for each $i\in F$. The \emph{gap function $g_{x^*}:F  \to \N$ of $x^*$}  is defined as follows:
\[ g_{x^*}(i) = \min I_{i+1} - \max I_i -1, \text{ if } i< \sup F, \text{ and } g_{x^*}(i)= 0 \text{ if } i = \sup F \]
\end {definition}

\begin{lemma}\label{lem5.7}
Let $x^*\in \mathcal{D}_1 \cap S_{J^*}$ such that $x^*= \sum_{i\in F}\alpha_i I_i^*$ with $\alpha_i\ne 0$ for all $i\in F$.
If $1 \in g_{x^*} [F]$ then $x^* \notin \overline{Ext(B_{J^*})}^{\|\cdot \|_*}$.
\end{lemma}

\begin{proof}
Since $\| x^*\|_* = 1$, Proposition \ref{prop5.1} yields that $x= \sum_{i \in F}\alpha_i e_i \in Ext(B_J)$. Therefore, for every $i < \sup F$ we have that  $\alpha_i \alpha_{i+1} < 0 $ (Proposition \ref{P02}).

Assume, on the contrary, that there is a sequence $(x_m^*)_{m\in \N} \subset Ext(B_{J^*})$ such that $x_m^* \xrightarrow{\|\cdot \|_*} x^* $. Choose $i_0 \in F$ such that $g_{x^*}( i_0) =1$ and set $n_0 = \max I_{i_0}$.
Then $n_0 +2= \min I_{i_0 +1}$ and, since $x_m^*  \xrightarrow{\| \cdot \|_*} x^*$, we conclude 
 \[ x_m^*(e_{n_0}) \to \alpha_{i_0}, \, x_m^*(e_{n_{0} +2}) \to \alpha_{i_{0} + 1} \text{ and } x_m^*(e_{n_{0} +1}) \to 0.\]

We set $x_m^*= \sum_{i \in F_m}\alpha_{i,m} I_{i,m}^* $. Theorem \ref{thm5.1} yields that $\cup _{i\in F_m} I_{i,m}$ is an interval and $x_m = \sum_ {i \in F_m} \alpha_{i,m} e_i \in Ext(B_J)$.
Thus, for $m\in \N$ large enough, $x_m^*(e_{n_0})  x_m^*(e_{n_{0} +2}) <0$  and either $x_m^*(e_{n_{0}+ 1}) = x_m^*(e_{n_0})$, or $x_m^*(e_{n_{0}+ 1}) =  x_m^*(e_{n_{0} +2})$, a contradiction since $x_m^*(e_{n_{0} +1}) \to 0$ and $\alpha_{i_0}, \alpha_{i_0+1} \ne 0$.
\end{proof}

Now we complete the description of $\overline{Ext(B_{J^*})}^{\| \cdot\|_*}$.

\begin{proposition}\label{prop5.8}
Let $x^*= \sum_{i \in F}\alpha_i I_i^* \in \mathcal{D}_1 \cap S_{J^*}$ such that $\alpha_i\ne 0$ for all $i \in F$. Then $x^*\in \overline{Ext(B_{J^*})}^{\|\cdot \|_*}$ if and only if $1 \notin g_{x^*}[F]$.
\end{proposition}
\begin{proof}
The direct implication is a consequence of Lemma \ref{lem5.7}.

For the converse, assume $1 \notin g_{x^*}[F]$, and set $G = \{ i\in F: g_{x^*}( i ) > 0 \}$. If $G = \emptyset$ then $x^* \in Ext(B_{J^*})$ (Proposition \ref{prop5.1}, Theorem \ref{thm5.1}). Thus we can assume that $G\ne \emptyset$.

We define a sequence $(x_m^*)_{m\in \N} \subset Ext(B_{J^*})$ such that $ x_m^*  \xrightarrow{\|_ \cdot \|_*} x^*$, and we proceed as in Example \ref{ex5.1}.
For $m\in \N$, we choose positive scalars $(\beta_i) _{i\in G}$ with $\sum_{i\in G} \beta_i^2 < 1/(2m^2)$, we decompose the interval $(\max I_i, \min I_{i+1})$ into two nonempty subintervals $L_i < M_i$ for each $i\in G$, and we define
$$
z_m^* = \sum_{i\in G} \bigl(-\sign (\alpha_i) \beta_iL_i^*+ \sign(\alpha_i)\beta_i M_i^*\bigr). 
$$

Since $\| z_m^* \|_* < 1/(\sqrt{2} m)$, setting $y_m^* = x^* + z_m^*$  and $x_m^* = y_m^*/\|y_m^*\|_*$, it is easy to check that $(x_m^*)_{m\in \N} \subset Ext(B_{J^*})$ and $\|x_m^*-x^*\|_*\to 0$.
\end{proof}

\section{Vectors with multiple $x$-norming partitions}\label{sect:multiple-part}

This is an appendix to Section \ref{sect:ext-B_J}. We provide examples of $x\in Ext(B_J)$ admitting multiple $x$-norming partitions, as it is stated in the following result.

\begin{proposition}\label{prop00}
Let $\alpha$ be either $k\in \N$ or $\aleph_0$ or $\mathfrak{c}$. Then there exists $x\in Ext(B_J)$ admitting exactly $\alpha$ $x$-norming partitions. These values of $\alpha$ are the only possible cardinalities of a set of $x$-norming partitions.
\end{proposition}

We will give the proof at the end of this section.

\begin{lemma}\label{lem0.1}
Let $x = a_1e_1 + a_2 e_2 + a_3e_3 \in \R^3$ such that $a_1> 0$, $a_2 < 0$ and $a_3 > 0$. Then the only possible $x$-norming partitions are $\mathcal{I} = \{1,2,3\}$ and $\mathcal{L} = \{ \{ 1 \}, \{2\},\{3\} \}$.
\end{lemma}
\begin{proof}
Since $a_1a_2<0$, $\{1,2\}$ is not an interval of an $x$-norming partition and the same happens for $\{a_2, a_3\}$.
\end{proof}

Let us denote
$$
E = \left\{ \sum_{i=1}^3 a_ie_i \in \R^3: \, a_1, a_3 > 0,\, a_2 <  0,   \, \sum_{i=1}^3 a_i = 1 \text{ and } \,  \sum_{i=1} ^3 a_i^2 =1\right\}.
$$
\medskip

\begin{lemma}\label{lem0.2}  The set $E$ satisfies the following properties:
\begin{enumerate}
\item $E$ is uncountable.
\item For every $0 < a_1 < 1$ there exists a unique pair $a_2, a_3$ in $\R$ such that
  \[ a_1e_1 + a_2e_2 + a_3e_3 \in E. \]
\item For every $a_1e_1 + a_2e_2 + a_3e_3 \in E$ we have $|a_2| < \min \{a_1, a_3\}$. If additionally $a_1 \ge 1/\sqrt{2}$ then $a_1 > a_3$.
  \end{enumerate}
\end{lemma}
\begin{proof}
(1) We consider the hyperplane $H = \{a_1e_1 + a_2 e_2 + a_3e_3  \in \R^3 : a_1 + a_2 + a_3 =1\}$ and we set $C = H \cap B_{\R^3} $. Then $C$ is the circumcircle of the triangle  with vertices $\{ e_1, e_2, e_3\} $.
We claim that $E$ is the set of points of the circumference of $C$ lying strictly between $e_1$ and $e_3$.

Indeed, the center of $C$ is $x_0= \sum_{i=1}^3 e_i/3$, and for $x\in C$ lying strictly between $e_1$ and $e_3$, the segments
$[x_0, x]$ and $[e_1, e_3] $ satisfy $[x_0, x] \cap [e_1, e_3] = \{ z\}$. Then
 \[z= sx_0 + (1-s)x = te_1 + (1-t)e_3,\]
hence $z(2) = s/3 + (1-s)x(2) = 0$. This yields $x(2) < 0$.
Moreover $x(1)$ and $x(3)$ are positive because $x(1)+x(2)+x(3) = 1 \, \text{ and }\, \| x\|_\infty \le 1.$

The same arguments show that no other vector in $C$ belongs to $E$.

(2) Observe that $E$ is connected and, if $P_1$ is the projection onto the first coordinate, then $P_1(E) = (0, 1)$, showing the existence of $a_2, a_3$. Moreover the hyperplane $P_1^{-1}(a_1)$ meets $E$ in exactly one point. This fact implies the uniqueness of the pair  $a_2, a_3$.

(3) It follows from the properties of $E$. Namely, $\sum_{i=1}^3 a_i = 1$ and $\sum_{i=1} ^3 a_i^2 =1$.
\end{proof}

Clearly, if $x\in E$ then $x\in Ext(B_J)$ and admits two different $x$-norming partitions.

\begin{lemma}\label{lem0.3} Let $G$ be an initial interval of $\N$ and let $(a_i)_{i\in G} \subset \R$ satisfying
  \[  \quad  a_1 > 0,\quad \sum_{i \in G\setminus \{1\}}a_i > 0 \quad \text{and}
    \quad |\sum_{i\in G} a_i|^2 = \sum_{i\in G} a_i^2. \]
Then for $ 0 < b_1 < a_1$ we have $| b_1 + \sum_{i \in G\setminus \{ 1\}} a_i|^2 < b_1^2 +\sum_{i\in G\setminus \{1\}} a_i^2 $.
\end{lemma}
\begin{proof}
From $|\sum_{i\in G} a_i|^2 = \sum_{i\in G} a_i^2$ we obtain
    \[a_1\bigl(\sum_{i\in G\setminus\{1\}} a_i\bigr) + \sum_{\substack{ i, j \ne 1, i \ne j } } a_i a_j =0\, \text{ and }\, a_1\bigl(\sum_{i\in G\setminus \{1\}} a_i \bigr)> 0.\]

These facts and $0 < b_1 < a_1$ yield
    \[ b_1\bigl(\sum_{i\in G\setminus\{1\}} a_i\bigr) + \sum_{\substack{ i, j \ne 1, i \ne j } } a_i a_j < 0, \]
which implies $ | b_1 + \sum_{i\in G\setminus\{1\}} a_i|^2 < b_1^2 +\sum_{i\in G\setminus\{1\}} a_i^2.$
\end{proof}

Our next result is the key to construct vectors $x\in J$ with multiple $x$-norming partitions.

\begin{proposition}\label{prop0.1}
$(a)$ Let $(r_i)_{i=1}^k$ be a finite strictly increasing sequence with $r_1 \ge 1/\sqrt{2}$ and $r_k = 1$.
Then there exists a unique finite sequence $( a_j)_{j=1}^{2k+1}$ satisfying
\begin{enumerate}
\item $a_1= r_1$, and for $1< m <k$,
    \[ \sum_{j=1}^{2m+1} a_j= r_{m +1}\, \text{ and }\, \bigl|\sum_{j=1}^{2m+1} a_j\bigr|^2 = \sum_{j=1}^{2m+1} a_j^2.\]
\item For every subinterval $I$ of $\{1,\ldots, 2k+1\}$ with $\#I >1$ and either $1< \min I$ or $\min I = 1$ and $\max I = 2n$, we have
     $|\sum_{j\in I} a_j|^2 <\sum_{j\in I} a_j^2 $.
\end{enumerate}
$(b)$ If $(r_k)_{k\in \N}$ is strictly increasing with $r_1 \ge 1/\sqrt{2}$ and $\lim r_k = 1$ then there exists a sequence
$(a_j)_{j\in \N}$ satisfying the analogous to the previous conditions (1), (2) in $(a)$.
\end{proposition}
\begin{proof}
$(a)$ We define $(a_j)_{j=1}^{2k+1}$ as follows. For $j=1$ we set $a_1= r_1$. For $m= 2,\ldots, k$, by part $(2)$ of Lemma \ref{lem0.2}, there are unique numbers $a_{2m}, a_{2m+1}$ satisfying
  \[ r_{m} + a_{2m} +  a_{2m+1} = r_{m+1}  \quad \text{and} \quad  r_{m}^2 + a_{2m}^2 +  a_{2m+1}^2 = r_{m+1}^2. \]
This completes the definition of  $(a_j)_{j=1}^{2k+1}$.
It is easy to check that
  \[ \sum_{j=1}^{2m+1} a_j = r_{m+1}\, \text{ and }\, |\sum_{j=1}^{2m+1 } a_j|^2 = \sum_{j=1}^{2m+1} a_j^2,\]
and part $(2)$ of Lemma \ref{lem0.2} implies that the sequence  $(a_j)_{j=1}^{2k+1}$ is unique.  So the proof of $(1)$ is complete.

To prove $(2)$, note that for $1\le m\le k-1$ the vector $ x_m= r_{m} e_1 +a_{2m} e_2+ a_{2m+1}e_3$ belongs to the set
  \[ E_m = \bigl\{\sum_{i=1}^3 a_ie_i \in \R^3: \, a_1, a_3 > 0,\, a_2 <  0,  \, \sum_{i=1}^3 a_i = r_{m+1} \text{ and }
             \sum_{i=1} ^3 a_i^2 =r_{m+1}^2\bigr\}.\]
Hence $a_{2m} < 0$ and $a_{2m+1} > 0$. Since $r_m \ge 1/\sqrt{2}$, fact (3) in Lemma \ref{lem0.2} yields
    \[ 0 < a_{2m+1} < r_m , \quad a_{2m} < 0 \quad \text{and} \quad |a_{2m}| < a_{2m+1}. \]

Let $I$ be a subinterval of $\{1, \ldots, 2k+1\}$ with $\#I >1$. First we assume that $1< 2m +1= \min I$ and $\max I=2l+1$. Then
$$
\sum_{j = 2m+1}^{2l+1} a_j=\bigl( \sum_{j = 1}^ {2l+1} a_j  -\sum_{j=1}^{2m+1} a_j \bigr) + a_{2m+1} = ( r_{l+1} - r_{m+1} )+ a_{2m+1} > 0
$$
and $\sum_{j=2m+2}^{2l+1} a_j= r_{l +1}- r_{m+1}> 0$. Therefore
$$
|r_{m+1}+ \sum_{j=2m+2}^{2l+1} a_j|^2= r_{m+1}^2 + \sum_{j=2m+2}^{2l+1} a_j^2 \quad \textrm{and} \quad r_{m+1} > a_{2m+1}.
$$

The last inequality holds since $r_{m+1} > r_m > a_{2m+1}$. Thus the assumptions of Lemma \ref{lem0.3} are fulfilled, and we get
 $|\sum _{j\in I} a_j |^2 < \sum_{j\in I} a_j^2$.
This completes the proof for subintervals of $\{1, \ldots, 2k-1\}$ with odd ends.

The other cases can be reduced to this one.
Indeed, assume first that $\min I = 2m+1$ and $\max I = 2l $. Then from the previous case we get $\sum_{j=2m+1}^{2l-1} a_j > 0$ and $a_{2l} < 0 $.
Hence
$$
|\sum_{j= 2m+1} ^{2l}a_j |^2 < |\sum _{j=2m+1}^{2l - 1} a_j|^2 + a_{2l}^2 \le  \sum _{j=2m+1}^{2l} a_j^2.
$$

With the same arguments we treat the case $\min I = 2m $ and $\max I = 2l+1$. The remaining case is  $\min I = 2m $ and $\max I = 2l$.

Again, as in the previous cases, we have $\sum_{j= 2m}^{2l-1}a_j > 0$ and $a_{2l} < 0$.
And similarly we conclude that
$$
|\sum_{j= 2m} ^{2l}a_j |^2 < |\sum _{j=2m}^{2l-1} a_j|^2 + a_{2l}^2 < \sum _{j=2m}^{2l} a_j^2.
$$

To finish the proof of $(2)$ of part $(a)$ we examine the case  $\min I=1$ and $\max I = 2n$, in which we have
$$
\sum_{j=1}^{2n}a_j = \sum_{j=1}^{2n-1}a_j + a_{2n}  = r_{n} + a_{2n}
$$
with $r_n>0$, $a_{2n}< 0$ and $|a_{2n}|< r_{n}$. Hence
$$
\sum_{j=1}^{2n}a_j |^2= |r_{n} + a_{2n}|^2 < |r_{n}|^2 <  \sum_{j=1}^{2n}a_j^2.
$$

To prove $(b)$, let $(r_k)_{k\in\N}$ as in the statement. We define the sequence $(a_j)_{j\in \N}$ in the same way as in the finite case, and satisfying
$$
\sum_{j=1}^{2k+1} a_j = r_{k+1}\, \text{ and }\, |\sum_{j=1}^{2k+1} a_j |^2 = \sum_{j=1}^{2k+1} a_j ^2  = r_{k+1}^2\, \text{ for all }\, k\in \N.
$$

Clearly  $a_j \to 0$. Hence  $\sum_{j=1}^\infty a_j = 1$ and $|\sum_{j=1}^\infty a_j |^2 = \sum_{j=1}^\infty a_j^2 = 1$.

The uniqueness of  $(a_j)_{j\in \N}$ can be proved as in the finite case. Thus part $(1)$ of $(b)$ is proved.
To prove $(2)$, we shall show that for every interval $I$ of $\N$ with $\#I >1$ and either $\min I > 1$ or $\min I = 1$ and  $\max I =2n$, we have
$| \sum _{j\in I} a_j |^2 < \sum_{j\in I} a_j^2$.

If $I$ is finite and satisfies one of the above two conditions then the proof is identical to the corresponding one in part $(a)$. Thus we assume that $I$ is infinite and $\min I > 1$.

If $\min I = 2m+1$, then
  \[ \sum_{j=2m+1}^\infty a_j = \sum_{j=2m+2}^\infty a_j + a_{2m+1} = 1-r_{m+1} + a_{2m+1} > 0\, \text{ and } \]

    \[ \sum_{j=2m}^\infty a_j= 1-r_{m}> 0\, \text{ and }\, r_{m} >  a_{2m+1}.\]
Hence Lemma \ref{lem0.3} yields the result.

If $\min I =2m$ then $a_{2m} < 0$ and $\sum_{j=2m+1}^\infty a_j> 0$. Hence
    \[ |  \sum_{j=2m}^\infty a_j |^2  < a_{2m} ^2 +| \sum_{j=2m+1}^\infty a_j |^2 <   \sum_{j=2m}^\infty a_j ^2, \]
part $(2)$ in $(b)$ is proved and the proof is complete.
\end{proof}

For $x\in J$, we denote $\mathcal{N}(x)=\{ \mathcal{I}:  \mathcal{I}\, \text{ is an $x$-norming partition}\}$.

\begin{corollary}\label{cor0.1} $(a)$ For every $k\in \N$ there exists $x_k \in Ext(B_J)$ such that $\#\mathcal{N}(x_k)=k$.

$(b)$ There exists $x_\infty \in Ext(B_J)$  such that $\#\mathcal{N}(x_\infty)= \aleph_0.$
\end{corollary}
\begin{proof}
(a) Let $k \in \N$. If $k=1$ then every $e_n$ satisfies the conclusion. If $k > 1$ we choose $(r_i)_{i=1}^k$ and set
$ x_k = \sum_{j=1}^{2k-1} a_j e_j$,
where $(a_j) _{j=1}^{2k-1} $ is the sequence associated to $(r_i)_{i=1}^k$ in Proposition \ref{prop0.1}.

Let $\{I_p \}_{p\in G}$ be an $x_k$-norming partition. Then Proposition \ref{prop0.1} yields that there is at most
one $I_p$ with $\#I_p > 1$ and, if it exists, then $\min I_p=1, \max I_p=2l-1$, $l=1,\ldots k$. Therefore $x_k \in Ext(B_J)$ and $\#\mathcal{N}(x_k)=k$.

$(b)$  The vector  $x_\infty $ is defined in a similar way.
\end{proof}

\begin{proposition}\label{prop0.2}
The set of all $x\in Ext(B_J)$ for which $\#\mathcal{N}(x) > \aleph_0$ is uncountable.
\end{proposition}
\begin{proof}
We take $x= a_1e_1 + a_2 e_2 + a_3e_3 \in E$ and we set $y = -a_3 e_4 - a_2  e_5  -a_1 e_6$.
Then we decompose $\N$ into disjoint intervals $(F_n)_{n\in \N}$ with $\#F_n = 6$ for all $n\in \N$.
We also denote $x_n + y_n$ the shift of $x+y$ into the interval $F_n$, and
    $z= \sum_{n=1}^\infty 2^{-n}(x_n+y_n). $

Let $\{I_i\}_{i\in G}$ be a $z$-norming partition. Then
\begin{enumerate}
\item For every $i\in G$ there exists $n\in \N$  such that
    \[\text{either }I_i \subset \supp(x_n) \, \text{ or } \, I_i \subset \supp(y_n).\]
\item For every $i \in G$ there exists $n\in \N$ such that
    \[\text{either } I_i= \supp(x_n)\, \text{ or }\, I_i = \supp(y_n)\, \text{ or }\, \exists k\in F_n\, \text{ with }\, I_i = \{k\}.\]
\end{enumerate}

Both assertions are easily established. Therefore, if we denote
      \[ \mathcal {I}_{x_n} = \{ \{k\} : k\in \supp(x_n) \}, \quad \mathcal {L}_{x_n} = \supp( x_n)\]
and we define
   \[ P = \prod_{n=1}^\infty \bigl \{ \mathcal {I}_{x_n}, \mathcal {L}_{x_n}\bigr\} \times\bigl\{ \mathcal{I}_{y_n}, \mathcal{L}_{y_n} \bigr\},\]
then every $\mathcal{I} \in P$ defines a $z$-norming partition, showing that $(\|z\|_J)^{-1}z\in Ext(B_J)$ and admits uncountably many
$z$-norming partitions. Since $z$ is different  for different vectors in $E$, the result is proved.
\end{proof}

\begin{remark} \label{rem0.1}
If for some $x\in J$ the set of $x$-norming partitions is uncountable, then this set has the cardinality of the continuum, as we showed in Remark \ref{rem3.1}.
%
\end{remark}

\begin{proof}[Proof of Proposition \ref{prop00}]
It is a direct consequence of Corollary \ref{cor0.1}, Proposition \ref{prop0.2} and Remark \ref{rem0.1}.
\end{proof}
\medskip

\noindent\textbf{Acknowledgments.}
We thank I. Gasparis for his help in this  research. In particular, his suggestion to use Choquet’s theorem in Section 5
made it more clear.
We also thank W.B. Johnson and P. Motakis for their comments. Finally we extend our thanks to the referee for his comments
and remarks enabling us to improve the presentation of the results. 

The second author was supported in part by MICINN (Spain), Grant PID2019-103961GB-C22.

\end{document}